\documentclass[12pt]{article}
\usepackage{amssymb,amsmath}
\author{Abdallah Assi\thanks{Universit\'e d'Angers, Math\'ematiques, 49045 Angers ceded 01,
France, e-mail:assi@univ-angers.fr\break Visiting address: American
University of Beirut, Department of Mathematics, Beirut 1107 2020,
Lebanon}}

\title{The tree model of a meromorphic plane curve}
\date{\mbox{}}

\newtheorem{teorema}{Theorem}[section]

\newtheorem{proposicion}[teorema]{Proposition}
\newtheorem{lema}[teorema]{Lemma}
\newtheorem{definicion}[teorema]{Definition}

\newtheorem{corolario}[teorema]{Corollary}

\newtheorem{nota}[teorema]{Remark}
\newtheorem{exemple}[teorema]{Example}

\newenvironment{demostracion}{\noindent Proof}

\newcommand{\RR}{{\mathbb R}}
\newcommand{\KK}{{\mathbb K}}

\begin{document}
\maketitle

\noindent{\bf Abstract.} We associate with a plane meromorphic
curve $f$ a tree model $T(f)$ based on its contact structure. Then
we give a  description of the $y$-derivative of $f$ (resp. the
Jacobian $J(f,g)$) in terms of $T(f)$ (resp. $T(fg)$). We also
characterize the regularity of $f$ in terms of its tree.

\section*{Introduction}

\medskip

\noindent  Let $\KK$ be an algebraically closed field of
characteristic $0$, and let $f,g$ be two monic reduced polynomial of
$\KK((x))[y]$ of degrees $n,m$. Let $f_x,g_x$
(resp. $f_y,g_y$) be the $x$-derivative (resp. the $y$-derivative) of
$f,g$, and let $J(f,g)=f_xg_y-f_yg_x$. Let, by Newton Theorem,

$$
f(x,y)=\prod_{i=1}^n(y-y_i(x)),\quad g(x,y)=\prod_{j=1}^m(y-z_j(x))
$$

\noindent where $(y_i(x))_{1\leq i\leq n}$ and $(z_j(x))_{1\leq j\leq
  m}$ are meromorphic fractional series in $x$.

\noindent The main objective of this paper is to give a description
of $f_y$ (resp. $J(f,g)$) when the contact structure of $f$ (resp.
$fg$) is given. Let $H(x,y)=\prod_{i=1}^a(y-Y_i(x))$ and
$\bar{H}(x,y)=\prod_{j=1}^b(y-Z_j(x))$ be two irreducible
polynomials of ${\KK}((x))[y]$ and define the contact $c(H,\bar{H})$
of $H$ with $\bar{H}$ to be

$$
c(H,\bar{H})={\rm max}_{i,j}O_x(Y_i-Z_j)
$$

\noindent  where $O_x$ denotes the $x$-order (in particular,
$c(H,H)=+\infty$). Let $f$ be as above and define the contact set of
$f$ to be

$$
C(f)=\lbrace O_x(y_i-y_j)|1\leq i\not= j\leq n\rbrace
$$

\noindent Let $f=f_1.\ldots.f_{\xi(f)}$ be
 the factorization of $f$ into irreducible components in
 $\KK((x))[y]$.
 Given $M\in C(f)$, we define $C_M(f)$ to be the set of irreducible components of $f$ such that $f_i\in C_M(f)$ if
 and only if $c(f_i,f_j)\geq M$ for some $j$ (with the understanding that $c(f_i,f_i)\geq M$ if
 and only if $M\geq O_x(y-y')$ for some roots $y\not=y'$ of $f_i(x,y)=0$). Given $f_i,f_j\in C_M(f)$, we say that $f_iR_Mf_j$
 if and only if $c(f_i,f_j)\geq M$. This defines
 an equivalence relation in $C_M(f)$.  The set of points of the tree of $f$ at the level $M$
is defined to be the set of equivalence classes of $R_M$. The set of
points defined this way -where two close points are connected with a
segment of line and top points are assigned with arrows- defines the
tree $T(f)$ of $f$:

\vskip1cm
 \unitlength=1cm
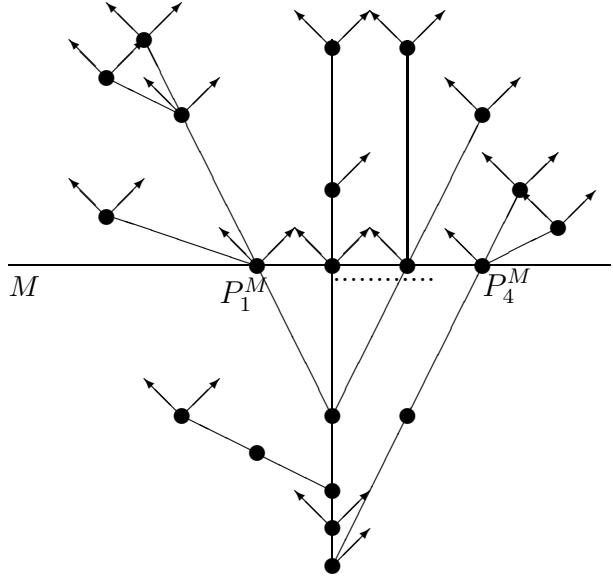
 \begin{center}
\begin{picture}(6,4)(3,-1)

\put(6,-4){\line(0,10){7}} 
\put(6,-4){\circle*{.2}} \put(6,-4){\vector(1,1){0.5}}
\put(6,-3.5){\circle*{.2}}
\put(6,-3.5){\vector(1,1){0.5}}\put(6,-3.5){\vector(-1,1){0.5}}
\put(6,0){\circle*{.2}}
\put(6,0){\vector(1,1){0.5}}\put(6,0){\vector(-1,1){0.5}}
\put(6,1){\circle*{.2}}\put(6,1){\vector(1,1){0.5}}
\put(6,2.9){\circle*{.2}}\put(6,2.9){\vector(1,1){0.5}}\put(6,2.9){\vector(-1,1){0.5}}

\put(6,-4){\line(1,2){2.5}}  
\put(7,-2){\circle*{.2}}
\put(8.5,1){\circle*{.2}}\put(8.5,1){\vector(1,1){0.5}}\put(8.5,1){\vector(-1,1){0.5}}

\put(6,-3){\line(-2,1){2}} 
\put(6,-3){\circle*{.2}} \put(5,-2.5){\circle*{.2}}
\put(4,-2){\circle*{.2}}\put(4,-2){\vector(1,1){0.5}}\put(4,-2){\vector(-1,1){0.5}}

\put(6,-2){\line(1,2){2}}  
\put(8,2){\circle*{.2}}
\put(8,2){\vector(1,1){0.5}}\put(8,2){\vector(-1,1){0.5}}

\put(6,-2){\circle*{.2}}
\put(6,-2){\line(-1,2){2.5}} 
\put(5,0){\circle*{.2}} \put(5,0){\line(-3,1){2}}
\put(5,0){\vector(-1,1){0.5}} \put(5,0){\vector(1,1){0.5}}
\put(3,0.65){\circle*{.2}}\put(3,0.65){\vector(1,1){0.5}}\put(3,0.65){\vector(-1,1){0.5}}

\put(4,2){\circle*{.2}}\put(4,2){\vector(1,1){0.5}}\put(4,2){\vector(-1,1){0.5}}

\put(3.5,3){\circle*{.2}}
\put(3.5,3){\vector(1,1){0.5}}\put(3.5,3){\vector(-1,1){0.5}}

 \put(4,2){\line(-2,1){1}} 
\put(3,2.5){\circle*{.2}}
\put(3,2.5){\vector(1,1){0.5}}\put(3,2.5){\vector(-1,1){0.5}}

\put(7,0){\line(0,3){3}} 
\put(7,0){\circle*{.2}} \put(7,0){\vector(-1,1){0.5}}
\put(7,2.9){\circle*{.2}}\put(7,2.9){\vector(1,1){0.5}}\put(7,2.9){\vector(-1,1){0.5}}

\put(8,0){\line(2,1){1}} \put(8,0){\circle*{.2}}
\put(8,0){\vector(-1,1){0.5}} \put(9,.5){\circle*{.2}}
\put(9,.5){\vector(1,1){0.5}}\put(9,.5){\vector(-1,1){0.5}}

\put(1.7,0){\line(1,0){8}} \put(1.7,-0.4){\shortstack{\small{$M$}}}
 \put(4.5,-0.5){\shortstack{$\small
P_1^M$}} \put(6,-0.2){\shortstack{$............$}}
\put(8,-0.4){$\displaystyle{\small{P_4^M}}$}
\end{picture}
\end{center}
\bigskip
\vskip3cm

\noindent Let $P_i^M$ be a point of the tree of $f$ at the level
$M$, and let $\bar{f}$ be a monic  polynomial of $\KK((x))[y]$. We
denote by $Q_{\bar{f}}(M,i)$ the product of irreducible components
of $\bar{f}$ whose contact
 with any element of $P_i^M$ is $M$. It results from [8] that deg$_yQ_{f_y}(M,i)>1$, i.e. every point
 of $T(f)$ gives rise to a component of $f_y$.   We give in Section 7.,
based on the results of Section 5., the $y$-degree of $Q_{f_y}(M,i)$
(see Proposition 7.6.), its intersection multiplicity  as well as
the contact of its irreducible components with $f_j,1\leq j \leq
\xi(f)$
 (see Theorem 7.7. and Theorem 8.9.). This result gives a generalization of Merle Theorem ($f\in\KK[[x]][y]$ and $\xi(f)=1$)
(see Proposition 7.1.) and Delgado Theorem ($f\in\KK[[x]][y]$ and
$\xi(f)=2$) (see Example 7.11.).
 These two results use the arithmetic
of the semigroup associated with $f$, which does not help for
meromorphic curves and, as shown by Delgado, does not seem to
suffice when $f\in\KK[[x]][y]$ and $\xi(f)\geq 3$.

\medskip

\noindent Let $T(fg)$ be the tree of $fg$. A point $P_i^M$ of
$T(fg)$ is said to be an $f$-point (resp. a $g$-point) if $P_i^M$
does not contain irreducible components of $g$ (resp. $f$). A point
of $T(fg)$ which is neither an $f$-point nor a $g$-point is called a
mixed point. This gives us the following description of $T(fg)$:

\vskip0.5cm
 \unitlength=1cm
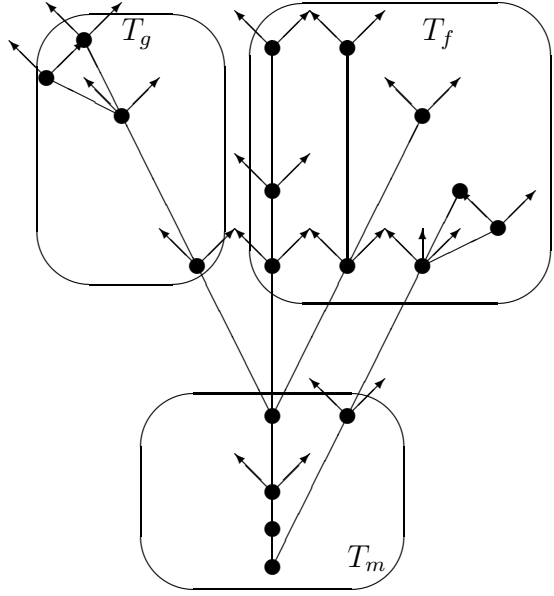
 \begin{center}
\begin{picture}(6,4)(3,-1)

\put(6,-4){\line(0,10){7}} 
\put(6,-4){\circle*{.2}} \put(7,-4){\shortstack{$T_m$}}
 \put(6,-3.5){\circle*{.2}}
\put(6,0){\circle*{.2}}\put(6,0){\vector(1,1){0.5}}\put(6,0){\vector(-1,1){0.5}}
\put(6,1){\circle*{.2}}\put(6,1){\vector(1,1){0.5}}\put(6,1){\vector(-1,1){0.5}}
\put(6,2.9){\circle*{.2}}\put(6,2.9){\vector(1,1){0.5}}\put(6,2.9){\vector(-1,1){0.5}}
\put(6,-4){\line(1,2){2.5}}  
\put(7,-2){\circle*{.2}}\put(7,-2){\vector(1,1){0.5}}\put(7,-2){\vector(-1,1){0.5}}
\put(8.5,1){\circle*{.2}}
\put(6,-3){\circle*{.2}}\put(6,-3){\vector(1,1){0.5}}\put(6,-3){\vector(-1,1){0.5}}
 \put(6,-3){\oval(3.5,2.6)}
\put(6,-2){\line(1,2){2}}  
\put(8,2){\circle*{.2}}\put(8,2){\vector(1,1){0.5}}\put(8,2){\vector(-1,1){0.5}}
\put(8,3){\shortstack{$T_f$}} \put(6,-2){\circle*{.2}}
\put(6,-2){\line(-1,2){2.5}} 
\put(5,0){\circle*{.2}}\put(5,0){\vector(1,1){0.5}}\put(5,0){\vector(-1,1){0.5}}
 \put(4,2){\circle*{.2}}\put(4,2){\vector(1,1){0.5}}\put(4,2){\vector(-1,1){0.5}}
\put(3.5,3){\circle*{.2}}\put(3.5,3){\vector(1,1){0.5}}\put(3.5,3){\vector(-1,1){0.5}}
\put(4,3){\shortstack{$T_g$}}
 \put(4,2){\line(-2,1){1}} 
\put(3,2.5){\circle*{.2}}\put(3,2.5){\vector(1,1){0.5}}\put(3,2.5){\vector(-1,1){0.5}}
\put(7,0){\line(0,3){3}} 
\put(7,0){\circle*{.2}}\put(7,0){\vector(1,1){0.5}}\put(7,0){\vector(-1,1){0.5}}
 \put(7.7,1.5){\oval(4,4)}
 \put(4.12,1.55){\oval(2.5,3.6)}
 \put(7,2.9){\circle*{.2}}\put(7,2.9){\vector(1,1){0.5}}\put(7,2.9){\vector(-1,1){0.5}}
\put(8,0){\line(2,1){1}}
\put(8,0){\circle*{.2}}\put(8,0){\vector(1,1){0.5}}\put(8,0){\vector(-1,1){0.5}}\put(8,0){\vector(0,1){0.5}}
\put(9,.5){\circle*{.2}}\put(9,0.5){\vector(1,1){0.5}}\put(9,0.5){\vector(-1,1){0.5}}
\end{picture}
\end{center}
\bigskip
\bigskip
\bigskip
\bigskip
\bigskip
\bigskip
\bigskip
\bigskip
\bigskip
\bigskip

 \noindent In Section 8, based on the results
of Sections 4. and 5., we prove the following:

\noindent{\bf Theorem} If $P_i^M$ is an $f$-point (resp. a
$g$-point), then deg$_yQ_{J(f,g)}(M,i)>1$.

\medskip

\noindent We also give an explicit formula for
deg$_yQ_{J(f,g)}(M,i)$ and its intersection multiplicity as well as
the contact  of  its irreducible components with each of the
irreducible components of $fg$ (see Theorem 8.4.). 

\medskip

\noindent As a consequence of this result, if
$J(f,g)\in{\mathbb K}((x))$, then every point of $T(fg)$ 
is a mixed point.

\medskip

\noindent Our explicit formulas for degrees, contacts and intersection multiplicities 
are given in terms of the invariants associated with the tree models of $f,g$ and $fg$. They are obtained 
using the results of Section 5 and Section 6.  Although these results are technical, we think that such 
precise formulas would be of interest for the study of problems such as the Jacobian conjecture in the
plane.
 
\medskip

\noindent The problem of the factorization of $f_y$ and $J(f,g)$ has been
considered by several  authors, with a special attention to the analytical case.
Beside the results of Merle and Delgado, Garc\'{i}a Barroso (see  [7]) used the Eggers tree
in order to get a decomposition of the generic polar of an analytic reduced curve 
(see  [12] for the definition and the properties of the Eggers tree).  In [9] and [10], Maugendre
computed the set of Jacobian quotients of a germ $(h_1,h_2):({\mathbb C}^2,0)\longmapsto ({\mathbb C}^2,0)$
in terms of the minimal resolution of $h_1h_2$.

\medskip

\noindent Let the notations be as above, and assume that
$f,g\in\KK[x^{-1}][y]$. Let $F(x,y)=f(x^{-1},y)$ and
$G(x,y)=g(x^{-1},y)$. For all $\lambda\in \KK$, we denote by
$F_{\lambda}$ the polynomial $F-\lambda$. We say that the family
$(F_{\lambda})_{\lambda\in\KK}$ is regular if the rank of the
$\KK$-vector space ${\dfrac{\KK[x,y]}{(F_{\lambda},F_y)}}$, denoted
Int$(F-\lambda,F_y)$, does not depend on $\lambda\in\KK$. When
$(F_{\lambda})_{\lambda\in\KK}$ is not regular, there exists a
finite number $\lambda_1,\ldots,\lambda_s\in\KK$ such that
Int$(F-\lambda,F_y)>{\rm Int}(F-\lambda_i,F_y)$ for $\lambda$
generic and $1\leq i\leq s$. The set $\lbrace
\lambda_1,\ldots,\lambda_s\rbrace$ is called the set of irregular
values of $(F_{\lambda})_{\lambda\in\KK}$.

\medskip

\noindent The regularity of a family of affine curves is related to
many problems in affine geometry, in particular the plane Jacobian
problem. If $(F_\lambda)_{\lambda}$ is regular and smooth,
then $F$ is equivalent to a coordinate of $\KK^2$. If
$(F_\lambda)_{\lambda}$ is smooth with only one irregular value
$\lambda_1$, then $F-\lambda_1$ is reducible in $\KK[x,y]$ and one of
its irreducible components is equivalent to a coordinate of
$\KK^2$. In general, nothing is known when $(F_\lambda)_{\lambda}$ has
more than two irregular values (see [4] and references).

\medskip

\noindent  Suppose that $F$ is generic in the family
$(F_\lambda)_{\lambda}$. In particular,  the intersection
multiplicity of $f$ with any irreducible component of $f_y$ is less
than $0$. Let  $P_i^M$ be a point of $T(f)$. We say that $P_i^M$  is
a bad point if one of the irreducible components of $Q_{f_y}(M,i)$
has intersection multiplicity $0$ with $f$. Otherwise, $P_i^M$ is
said to be a good point. Hence the tree $T(f)$ can be partitioned
into bad and good points. In Section 9 we characterize the notion of
regularity in terms of this partition. This, with the results of
Section 2. is used in Section 10. in order to prove that the set of
irregular values of $f$ is bounded by the number of irreducible
components $\xi(f)$ of $f$ (or equivalently the set of irregular
values of $(F_{\lambda})_{\lambda\in\KK}$ is bounded by the number
of places of $F$ at infinity).

\medskip

\noindent The author would like to think the referees for their valuable comments
and suggestions.

\section {Characteristic sequences}

\medskip

\noindent In this Section we shall recall some well known results
about the theory of meromorphic curves (see [2] for example). Let

 $$
 f= y^n +
a_1(x)y^{n-1} + ... + a_n(x)
$$

\noindent   be a monic  irreducible  polynomial of $\KK((x))[y]$,
where $\KK((x))$ denotes the field of meromorphic series over $\KK$.
Let, by Newton Theorem, $y(t)\in\KK((t))$ such that $f(t^n,y(t))=0$.
If $w$ is a primitive $n$th root of unity, then we have:

$$
f(t^n,y)=\prod_{k=1}^n (y-y(w^kt)).
$$

\noindent Write $y(t)=\sum_ia_it^i$, and let supp$(y(t))=\lbrace i;
a_i\not=0\rbrace$. Clearly supp$(y(t))={\rm supp}(y(w^kt))$ for all
$1\leq k\leq n-1$. We denote this set by supp$(f)$ and we recall
that gcd$(n,{\rm supp}(f))=1$. If we write $x^{1 \over n}$ for $t$,
then $y(x^{1\over n})= \sum a_i x^{i\over n}$ and $f(x,y(x^{1\over
n}))=0$, i.e. $y(x^{1\over n})$ is a root of $f(x,y)=0$. By Newton
Theorem, there are $n$ distinct roots of $f(x,y)=0$, given by
$y(w^kx^{1\over n}), 1\leq k\leq n$. We denote the set of roots of
$f$ by Root$(f)$.

\noindent We shall associate with  $f$ its characteristic sequences
$(m^f_k)_{k\geq 0}, (d^f_k)_{k\geq 1}$ and $(r^f_k)_{k\geq 0}$
defined by:

$\mid m^f_0\mid =d^f_1=\mid r^f_0\mid =n$, $m^f_1=r^f_1=$
inf($\lbrace i\in{\rm supp}(f)|{\rm gcd}(i,n)<{\rm min}(i,n)\rbrace$, and for all $k\geq 2$,

$d^f_k=$ gcd $(m^f_0,\ldots,m^f_{k-1})=$ gcd $(d^f_{k-1},
m^f_{k-1})$,

$m^f_k=$ inf $\lbrace i\in{\rm supp}(f)|$ $i$ is not divisible by
$d^f_k \rbrace$,

\noindent and $r^f_k=r^f_{k-1} \displaystyle{{d^f_{k-1}\over
d^f_k}}+m^f_k-m^f_{k-1}$.

\noindent Since gcd($n$, supp($f$)) $=1$, then there is $h_f\in
\mathbb N$ such that $d_{h_f+1}=1$. We denote by convention
$m^f_{h_f+1}=r^f_{h_f+1}=+\infty$. The sequence $(m_k)_{0\leq k\leq h_f}$ is also called
the set of Newton-Puiseux exponents of $f$. We finally set
$\displaystyle{e^f_k={d^f_k\over d^f_{k+1}}}$ for all $1\leq k\leq
h_f$.

\noindent Let $H$ be a polynomial of $\KK((x))[y]$. We define the
intersection of $f$ with $H$, denoted int$(f,H)$, by
int$(f,H)=O_tH(t^n,y(t))=n.O_xH(x,y(x^{1\over n}))$, where $O_t$
(resp. $O_x$) denotes the order in $t$ (resp. in $x$).

\noindent Let $p,q\in \mathbb{N}^*$, and let
$\alpha(x)\in\KK((x^{1\over
    p})), \beta(x)\in\KK((x^{1\over q}))$. We set

$$
{\rm c}(\alpha,\beta)=O_x(\alpha(x)-\beta(x))
$$

\noindent and we call  ${\rm c}(\alpha,\beta)$ the contact of
$\alpha$ with $\beta$. We define the contact of $f$ with $\alpha(x)$ to be

$$
{\rm c}(f,\alpha)=\displaystyle{{\rm max}_{1\leq i\leq n}
O_x(y_i(x)-\alpha(x))}
$$

\noindent where $\lbrace y_1,\ldots,y_n\rbrace={\rm Root}(f)$.

\noindent Let $g=y^m+b_1(x)y^{m-1}+\ldots+b_m(x)$ be a monic
irreducible polynomial of $\KK((x))[y]$ and let ${\rm Root}(g)=\lbrace
z_1,\ldots,z_m\rbrace$. We define
the contact of $f$ with $g$ to be

$$
{\rm c}(f,g)={\rm c}(f,z_1(x)).
$$

\noindent Note that ${\rm c}(f,g)= {\rm c}(f,z_j(x))={\rm
c}(g,y_i(x))$ for all $1\leq j\leq m$ and for all $1\leq i\leq n$.

\medskip

\begin{nota}{\rm (see [1]) i) Let $f\in\KK[[x]][y]$
(resp. $f\in\KK[x^{-1}][y]$). The set of int$(f,g)$,
$g\in\KK[[x]][y]$ (resp. $g\in\KK[x^{-1}][y]$) is a subsemigroup of
$\mathbb Z$. We denote it by  $\Gamma (f)$ and we call it the semigroup
associated with $f$. With the notations
above, $r^f_k >0$ (resp. $r^f_k <0$) for all $k=0,\ldots,h_f$, and
$r^f_0, r^f_1,\ldots,r^f_{h_f}$ generate $\Gamma(f)$. We write
$\Gamma(f)=<r^f_0, r^f_1,\ldots,r^f_{h_f}>$.

\noindent ii) For all $1\leq k \leq h_f$, $e^f_k$ is the minimal
integer such that $e^f_{k}r^f_{k} \in <r^f_0,r^f_1,\ldots,
r^f_{k-1}>$.

\noindent iii) For all $1 \leq k \leq h_f$, there is a monic
irreducible  polynomial $g_k \in \KK((x))[y]$ of degree
$\displaystyle{n \over d^f_k}$ in $y$ such that
c$(f,g_k)=\displaystyle{m^f_k \over n}$ and {\rm
int}$(f,g_k)=r^f_k$. Furthermore, $\displaystyle{\Gamma(g_k)= <
{r^f_0\over d^f_k},{r^f_1\over d^f_k},\ldots,{r^f_{k-1}\over
d^f_k}>}$. }\end{nota}

\begin{lema} {\rm (see [1]) Let $y(x)=\sum_ia_ix^{i\over n}\in {\rm Root}(f)$.
Given $s\in {\mathbb N}^*$, let $U_s$ denotes the group of the $s$th
roots of unity in $\KK$. Set

$R(i)=\lbrace w\in U_n |c(y(x),y(wx))=O_x(y(x)-y(wx))\geq
\displaystyle{m^f_i\over n}\rbrace$

$S(i)=\lbrace w\in U_n | c(y(x),y(wx))=
O_x(y(x)-y(wx))=\displaystyle{m^f_i\over n}\rbrace$.

\noindent We have the following:

i) For all $1\leq i \leq h_f+1$, $R(i)=U_{d^f_i}$. In particular,
{\rm card}$(R(i))=d^f_i$.

ii) For all $1\leq i \leq h_f$,
$S(i)=R(i)-R(i+1)=U_{d^f_i}-U_{d^f_i+1}$. In particular, {\rm
card}$(S(i))=d^f_i-d^f_{i+1}$.}\end{lema}

\begin{demostracion}{.}  Let $w\in U_n$, then
$y(x)-y(wx)=\sum_ka_k(1-w^k)\displaystyle{x^{k\over n}}$. In particular,
$O_x(y(x)-y(wx))\geq \displaystyle{m^f_i\over n}$ if and only if
$w^k=1$ for all $k < m^f_i$. This holds if and only if
$w\in U_{{d^f_i}}$.$\blacksquare$\end{demostracion}
\medskip

\begin{nota}{\rm i) Let $F$ be a nonzero monic polynomial of
      $\KK((x))[y]$. Assume that $F$ is reduced and let
    $F=F_1.\ldots.F_{\xi(F)}$ be the factorization of $F$ into
    irreducible polynomials of $\KK((x))[y]$. We define Root$(F)$
    to be the union of Root$(F_i),i=1,\ldots,\xi(F)$. Given a
    polynomial $G\in\KK((x))[y]$, we set
    int$(F,G)=\sum_{i=1}^{\xi(F)}{\rm int}(F_i,G)$.

ii) Let $p\in \mathbb{N}^*$, and let $F$ be a nonzero monic
polynomial of $\KK((x^{1\over p}))[y]$.  Assume that $F$ is reduced
and let
      $x=X^p, y=Y$, and $\bar{F}(X,Y)=F(X^p,Y)$. The polynomial $\bar{F}$ is
      a monic reduced polynomial of $\KK((X))[Y]$. Let
    Root$(\bar{F})=\lbrace Y_1(X),\ldots,Y_N(X)\rbrace$. The set of
    roots of $F(x,y)=0$ is $\lbrace
    Y_1(x^{1\over p}),\ldots,
Y_N(X^{1\over p})\rbrace$. }
\end{nota}

\noindent Let $M$ be a given real number and consider the sequence
$(m^f_k)_{1\leq k\leq h_f+1}$ of Newton-Puiseux exponents of $f$.  We define
the function $S(m^f,M)$ by putting

\[
S(m^f,M) =
\begin{cases}
r^f_k d^f_k + (nM-m^f_k)d^f_{k+1} & \text{ if  ${m^f_1\over n}\leq {m^f_k\over n} \leq M < {m^f_{k+1}\over n}$} \\
Md_1  & \text{if $M < {m^f_1\over n}$}
\end{cases}
\]

\begin{proposicion}{\rm (see [1] or [8]) Let $g=y^m+b_1(x)y^{m-1}+\ldots+b_m(x)$ be
a monic irreducible  polynomial of $\KK((x))[y]$.  We have the
following:

$c(f,g)=\displaystyle{M}$ \quad  if and only if \quad  {\rm
int}$(f,g)= S(m^f,M)\displaystyle{m\over n}$

 $c(f,g) < \displaystyle{M}$ \quad  if and only if\quad  {\rm int}$(f,g)< S(m^f,M)\displaystyle{m\over n}$

$c(f,g)>\displaystyle{M}$ \quad  if and only if \quad  {\rm
int}$(f,g)>S(m^f,M)\displaystyle{m\over n}$

}
\end{proposicion}















\medskip

\noindent Let $g_1,g_2$ be two monic irreducible  polynomials of
$\KK((x))[y]$ of degrees $q_1$ and $q_2$ respectively and let
$(m^{g_i}_k)_{1\leq k\leq h_{g_i}}$ be the set of Newton-Puiseux
exponents  of $g_i$, $i=1,2$.

\medskip

\begin{lema}{\rm  (see [1]) Let $\displaystyle{M}= {\rm
      min}(c(f,g_2),c(f,g_1))$. We have the following:

\noindent (i) $c(g_1,g_2) \geq \displaystyle{M}$.

\noindent (ii) if $c(f,g_2)\not=c(f,g_1)$ then
$c(g_1,g_2)=\displaystyle{M}$.
}
\end{lema}





\begin{lema} {\rm Let the notations be as above and let $(m_k^g)_{1\leq
k\leq h_g+1}$ be the set of Newton-Puiseux exponents of
$g$. Let $M=c(f,g)$ and assume that   $M\geq \dfrac{m^f_1}{n}$.  Let  $k$ be the greatest integer such that $\dfrac{m^f_k}{n}=\dfrac{m^g_k}{m}\leq M$. We have the following:

i) $\dfrac{n}{d^f_i}=\dfrac{m}{d_i^g}$ for all $i=1,\ldots,k+1$.

ii)  $\dfrac{n}{d^f_{k+1}}$ divides $m$. In particular, if $k=h$ then $n$ divides $m$.
}\end{lema}

\begin{demostracion}{.}  ii) results from i), since by i), $m=\dfrac{n}{d^f_{k+1}}d^g_{k+1}$. On the other hand, let $1\leq i\leq k$ and remark that
$m.n=n.m,m.m^f_1=n.m^g_1,\ldots, m.m^f_{i-1}=n.m^g_{i-1}$, in particular
$m.d^f_i=m.{\rm gcd}(n,m^f_1,\ldots,m^f_{i-1})=n.{\rm gcd}(q,m^g_1,\ldots,m^g_{i-1})=n.d^g_i$. This proves i).$\blacksquare$
\end{demostracion}

\medskip

\begin{lema}{\rm  Let the notations be as in Lemma 1.6. and let $y(x)\in{\rm Root}(f)$ (resp. $z(x)\in{\rm Root}(g)$)  such that $c(y(x),z(x))=M$. Write $y(x)=\sum_{i}c^f_ix^{i\over n}$ and $z(x)=\sum_{j}c^g_jx^{j\over m}$.
If  $M=\dfrac{m^f_{h_f}}{n}$ and $n \geq m$, then either $c^g_{mM}$ -the coefficient of $x^{M}$ in $z(x)$- is $0$, or $m=n$.
}
\end{lema}

\begin{demostracion}{.}  If $c^g_{mM}\not=0$, then
  $M=\dfrac{m^g_{h_g}}{m}$, hence $n$
divides $m$. This, with the hypotheses implies that $m=n$.$\blacksquare$
\end{demostracion}

\noindent As a corollary we get the following:

\begin{lema}{\rm Let $g_1,g_2$ be two monic irreducible  polynomials of $\KK((x))[y]$ of
degrees $q_1, q_2$ respectively,
and assume that $c(g_1,f)=c(g_2,f)=\dfrac{m^f_{h(f)}}{n}$. If $q_1 <n$
and $q_2 < n$, then $c(g_1,g_2) > \dfrac{m^f_{h_f}}{n}$.
}
\end{lema}

\begin{demostracion}{.} Let $y(x)\in{\rm Root}(f)$ (resp. $z_1(x)\in{\rm Root}(g_1)$, $z_2(x)\in{\rm Root}(g_2)$)
such that $c(y(x),z_1(x))=c(y(x),z_2(x))=\dfrac{m^f_{h_f}}{n}$. In particular
$c(z_1(x),z_2(x)) \geq \dfrac{m^f_{h(f)}}{n}$. By Lemma 1.7., the coefficients
of $x^{\dfrac{m^f_{h_f}}{n}}$ in $z_1(x)$ and $z_2(x)$ are $0$, which implies that
 $c(z_1(x),z_2(x)) > \dfrac{m^f_{h_f}}{n}$. This proves our assertion.$\blacksquare$
\end{demostracion}


\section{Equivalent and almost equivalent polynomials}

\noindent Let $f,g$ be two monic irreducible polynomials of $\KK((x))[y]$,
of degrees $n,m$ in $y$. Let
$(m^f_k)_{1\leq k\leq h_f}$, $(d^f_k)_{1\leq k\leq h_f}$, and $(r^f_k)_{0\leq k\leq h_f}$
(resp. $(m^g_k)_{1\leq k\leq h_g}$, $(d^g_k)_{1\leq k\leq h_g}$, and $(r^g_k)_{0\leq k\leq h_g}$) be the set of characteristic
sequences of $f$ (resp. of $g$).


\begin{definicion}{\rm i) We say that $g$ is equivalent to $f$ if the following holds:

- $h_f=h_g$

- $\displaystyle{{m^g_k\over m}={m^f_k\over n}}$ for all
$k=1,\ldots,h_f$.

- c$(f,g)\geq \displaystyle{{m^f_{h_f}\over n}}$.

\noindent ii) We say that $g$ is almost equivalent to $f$ if the
following holds:

- $h_f=h_g+1$.

- $\displaystyle{{m^f_k\over n}={m^g_k\over m}}$ for all $k=1,\ldots,
h_g$.

- c$(f,g)=\displaystyle{{m^f_{h_f}\over n}}$. }\end{definicion}

\begin{lema}{\rm Let the notations be as in Definition 2.1.

i) If $g$ is equivalent to $f$, then $m=n$.

ii) If $g$ is almost equivalent to $f$, then
$m=\displaystyle{{n\over d^f_{h_f}}}$. Furthermore, if
$y(x)=\sum_pc_px^{p\over m}\in{\rm Root}(g)$, then
$c_{{m^f_{h_f}\over n}.m}=0$.
 }\end{lema}

\begin{demostracion}{.} i) results from Lemma 1.6. On the other hand, by the
same Lemma, $m=a\displaystyle{{{n}\over{d^f_{h_f}}}}$ for
some $a\in\mathbb{N}^*$, but gcd$(a\displaystyle{{{n}\over{d^f_{h_f}}}}, \displaystyle{a\over {d^f_{h_f}}}m^f_1,\ldots,\displaystyle{a\over {d^f_{h_f}}}m^f_{h_f-1})=
\displaystyle{a\over d^f_{h_f}}d^f_{h_f}=1$, hence $a=1$. This proves
the first assertion of ii). Now the least assertion results from Lemma 1.7.$\blacksquare$
\end{demostracion}

\begin{definicion}{\rm
Let $\lbrace F_1,\ldots, F_r\rbrace$ be a set of monic irreducible
polynomials of $\KK((x))[y]$. Assume that $r > 1$ and let
$n_{F_i}={\rm deg}_yF_i$ for all $1\leq i\leq r$.

\noindent i) We say that the sequence $(F_1,\ldots,F_r)$ is
equivalent if for all $1\leq i\leq r$, $F_i$ is equivalent  to $F_1$.

\noindent ii) We say that the sequence $(F_1,\ldots,F_r)$ is
almost equivalent  if the following holds:

- The sequence contains an equivalent subsequence of $r-1$ elements.

- The remaining element is almost equivalent  to the elements of
the subsequence. }\end{definicion}

\begin{proposicion}{\rm Let the notations be as in Definition 2.3. and let $M$ be
a rational number. If $c(F_i,F_j)=M$ for all $i\not= j$, then the sequence
$(F_1,\ldots,F_r)$ is either equivalent or almost equivalent.
}
\end{proposicion}

\begin{demostracion}{.} If $r=1$, then there is nothing to prove. Assume that $r>1$, and that
$n_{F_1}={\rm max}_{1\leq k\leq r}n_{F_k}$.

- If $M > m^{F_1}_{h_{F_1}}$, then, by Lemma 1.6., ii), $n_{F_1}$ divides $n_{F_k}$
for all $1\leq k\leq r$. In particular $n_{F_1}=n_{F_k}$ and $F_k$ is equivalent to $F_1$
for all $1\leq k\leq r$.

- Suppose that $M=\displaystyle{{m^{F_1}_{h_{F_1}}}\over
{n_{F_1}}}$, and that $(F_1,\ldots,F_r)$ is not equivalent. Suppose,
without loss of generality, that $F_2$ is not equivalent to $F_1$.
By hypothesis, $M\geq \displaystyle{{m^{F_2}_{h_{F_2}}}\over
{n_{F_2}}}$ and $\displaystyle{{{m^{F_1}_j}\over
{n_{F_1}}}={{m^{F_2}_j}\over {n_{F_2}}}}$ for all $1\leq j\leq
h_{F_1}-1$. Let $y(x)=\sum c_px^{p\over n_{F_2}}\in {\rm
Root}(F_2)$. If the coefficient of $x^M$ in $y(x)$ is non zero, then
$n_{F_1}$ divides $n_{F_2}$, in particular $n_{F_2}=n_{F_1}$, and
$\displaystyle{{{m^{F_1}_{h_{F_1}}}\over
{n_{F_1}}}={{m^{F_2}_{h_{{F_2}}}}\over {n_{F_2}}}}$. Hence $F_1$ is
equivalent to $F_2$, which is a contradiction. Finally
$h_{F_2}=h_{F_1}-1$, and $n_{F_2}=a.\displaystyle{n_{F_1}\over
d^{F_1}_{h_{F_1}}}$, but
gcd$(n_{F_2},m^{F_2}_1,\ldots,m^{F_2}_{h_{F_2}})=1$, hence $a=1$ and
$n_{F_2}=\displaystyle{{n_{F_1}}\over d^{F_1}_{h_{F_1}}}$.  In
particular $F_2$ is almost equivalent to $F_1$. Let $k>2$. If $F_k$
is not equivalent to $F_1$, then $n_{F_k}=n_{F_2}<n_{F_1}$ by the
same argument as above. In particular, by Lemma 1.8., $c(F_1,F_2)>
M$, which is a contradiction. Finally the sequence
$(F_1,\ldots,F_r)$ is almost equivalent.$\blacksquare$
\end{demostracion}

\section{The Newton polygon of a meromorphic plane curve}

\medskip

\noindent In this Section we shall recall the notion of the Newton
polygon of a meromorphic plane curve. More generally let $p\in
\mathbb{N}$ and let $F=y^N+A_1(x)y^{N-1}+\ldots+A_{N-1}(x)y+A_N(x)$
be a reduced polynomial of $\KK((x^{1/p}))[y]$. For all
$i=0,\ldots,N$, let $\alpha_i=O_xA_i(x)$. The Newton boundary of $F$
is defined to be the boundary of the convex hull of
$\bigcup_{i=1}^N(\alpha_i,i)+{\RR_+}$.

\noindent Write $F(x,y)=\sum_{ij}c_{ij}x^{{i}\over{p}}y^j$ and
let Supp$(F)=\lbrace (\dfrac{i}{p},j)| c_{ij}\not=0\rbrace$, then
the Newton boundary of $F$ is also the boundary of the convex hull
of $\bigcup_{({{i}\over {p}},j)\in{\rm
Supp}(F)}(\dfrac{i}{p},j)+{\RR_+}$.

\noindent We define the Newton polygon of $F$, denoted $N(F)$, to be the union of the compact faces of the Newton boudary of $F$. Let $\lbrace P_k=(\alpha_{k_j},k_j), k_0 > k_1\ldots
>k_{v_F}\rbrace$   be the set of vertices of $N(F)$.  We denote this
set by $V(F)$. We denote by
$E(F)=\lbrace \triangle^F_l= P_{k_{l-1}}P_{k_{l}},
l=1,\ldots,v_F\rbrace$ the set of edges of $N(F)$. For all $1\leq
l\leq v_F$ we set $F_{\triangle^F_l}=\sum_{({i\over p},j)\in {\rm Supp}(F)\bigcap \triangle^F_l}c_{ij}x^{i\over p}y^j.
$
\bigskip

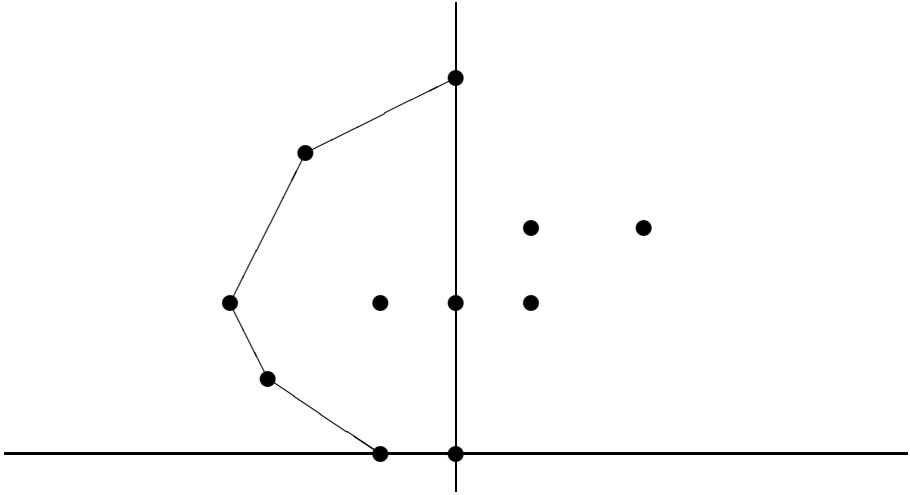
\unitlength=1cm
\begin{picture}(6,4)

\put(6,-2.5){\line(0,10){6.5}} 
\put(3,0){\circle*{.2}} \put(6,0){\circle*{.2}}
\put(8.5,1){\circle*{.2}}
\put(3,0){\line(1,-2){0.5}} \put(3.5,-1){\circle*{.2}}
\put(3.5,-1){\line(3,-2){1.5}} \put(5,-2){\circle*{.2}}
\put(3,0){\line(1,2){1}}
\put(6,-2){\circle*{.2}}
\put(5,0){\circle*{.2}} \put(4,2){\circle*{.2}}
\put(6,3){\circle*{.2}}
\put(4,2){\line(2,1){2}} 
\put(7,0){\circle*{.2}} \put(7,1){\circle*{.2}}
\put(0,-2){\line(1,0){12}}
\end{picture}
\bigskip
\bigskip
\bigskip
\bigskip
\bigskip
\bigskip
\bigskip

\begin{lema}{\rm Given $1\leq l\leq v_F$, there is exactly
    $k_{l-1}-k_l$ elements of Root$(F)$,
    $y_j(x), 1\leq j\leq k_{l-1}-k_l$, such that
     the following hold

i)
$O_x(y_j(x))=\displaystyle{\alpha_{k_{l-1}}-\alpha_{k_l}\over{k_{l-1}-k_l}}$
for all $j$.

ii) The set of initial
coefficients, denoted inco,  of $y_1,\ldots,y_{(k_{l-1}-k_l)}$ is nothing but the
set of nonzero roots of $F_{\triangle^F_l}(1,y)$.

\noindent Conversely, given $y(x)\in{\rm Root}(F)$, there exists $\triangle^F_l$ such that
$O_x(y(x))=\displaystyle{\alpha_{k_{l-1}}-\alpha_{k_l}\over{k_{l-1}-k_l}}$.

\noindent We denote the set of $x$-orders of ${\rm Root}(F)$ by $O(F)$, and we set Poly$(F)=\lbrace F_{\triangle^F_l}(1,y)|
1\leq l\leq v_F\rbrace$}.
\end{lema}

\begin{lema}{\rm Let $F$ be as above, and let $M$ be a rational
    number. Define $L_M:{\rm Supp}(F)\longmapsto \mathbb{Q}$
 by $L_M(\dfrac{i}{p},j)=\dfrac{i}{p}+Mj$, and let $a_0={\rm inf}(L_M({\rm
   Supp}(F)))$. Let
 in$_M(F)=\sum_{{i\over p}+Mj=a_0}c_{ij}x^{i\over p}y^j$. We have
 the following:

i) $M\in O(F)$ if and only if in$_M(F)$ is not a monomial. In this
case,
$M=\displaystyle{\alpha_{k_{l-1}}-\alpha_{k_l}\over{k_{l-1}-k_l}}$ for
some $1\leq l\leq v_F$, and in$_M(F)=F_{\triangle^F_l}$. Furthermore,
$(a_0,0)$ is the point where the line defined by
$(\alpha_{k_{l-1}},k_{l-1})$ and $(\alpha_{k_l},k_l)$ intersects the $x$-axis.

ii) Consider the change of variables $x=X,y=X^MY$ and let
$\bar{F}(X,Y)=F(X,X^MY)$. We have $\bar{F}=\sum c_{ij}x^{{i\over
    p}+Mj}=x^{a_0}F_{\triangle^F_l}(1,y)+\sum_{a > a_0}x^aP_a(y)$.
}
\end{lema}
\begin{demostracion}{.} Easy exercise.$\blacksquare$
\end{demostracion}


\noindent The following two lemmas give information about the Newton
polygons of the $y$-derivative (resp. the Jacobian) of a meromorphic
curve (resp. the Jacobian of two meromorphic curves).

\begin{lema}{\rm Let $F$ be as above and let $N(F)$ be the Newton
polygon of $F$. Let $V(F)=\lbrace P_k=(\alpha_{k_l},k_l), k_0 > k_1\ldots
>k_{v_F}\rbrace$ be the set of vertices of $F$ and assume that
$k_{v_F}=0$, i.e. $N(F)$
meets the $x$-axis. Assume that $(\alpha^1,1)\in {\rm
  Supp}(F_{\triangle^F_{v_F}})$ for some $\alpha^1\in\mathbb{Q}$, and
that $(\alpha^1,1)\notin V(F)$. We have the following

i) $(\alpha^1,0)\in V(F_y)$.

ii) $N(F_y)$ is the translation of $N(F)$
with respect to the vector $(0,-1)$.

iii) $O(F_y)=O(F)$, $v_F=v_{F_y}$.

iv)  deg$_yF_{\triangle^F_{v_F}}={\rm
    deg}_y({F_Y})_{\triangle^{F_y}_{v_F}}+1$. In particular, if $F$
  has $s$ roots whose order in $x$ is
  $\displaystyle{{\alpha_{k_{v(F)-1}}-\alpha_{k_{v(F)}}}\over{k_{v_F-1}}}$,
    then $F_y$ has $s-1$ roots with the same order in $x$.}
\end{lema}

\begin{demostracion}{.} The proof follows immediately from the
  hypotheses and Lemma 3.1.$\blacksquare$
\end{demostracion}

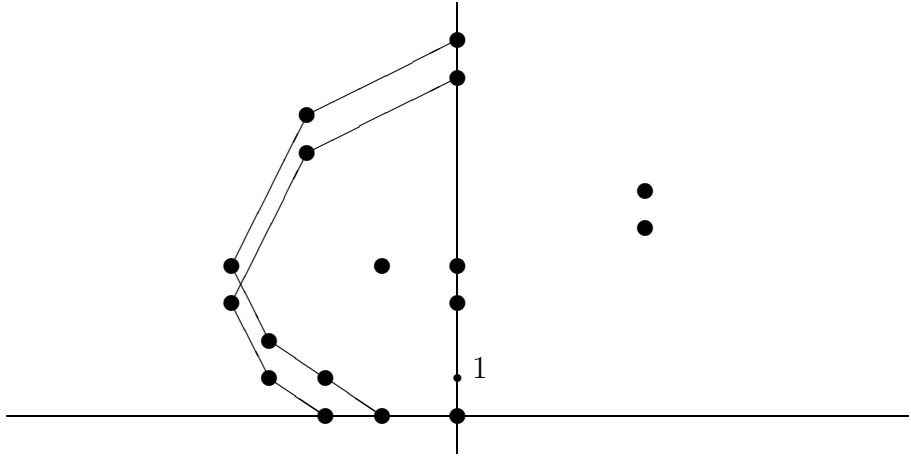
\unitlength=1cm
\begin{center}
\begin{picture}(6,4)(3,-0.5)

\put(6,-2.5){\line(0,10){6}} 
\put(3,0){\circle*{.2}} \put(6,0){\circle*{.2}}
\put(3,-0.5){\circle*{.2}} \put(6,-0.5){\circle*{.2}}
\put(8.5,1){\circle*{.2}}
\put(8.5,0.5){\circle*{.2}}
\put(3,0){\line(1,-2){0.5}}
 \put(3.5,-1){\circle*{.2}}
\put(3,-0.5){\line(1,-2){0.5}}
 \put(3.5,-1.5){\circle*{.2}}
\put(3.5,-1){\line(3,-2){1.5}}
\put(3.5,-1.5){\line(3,-2){0.75}}
 \put(5,-2){\circle*{.2}}
 \put(4.25,-1.5){\circle*{.2}}
 \put(4.25,-2){\circle*{.2}}
\put(3,0){\line(1,2){1}}
\put(3,-0.5){\line(1,2){1}}
\put(6,-2){\circle*{.2}}
\put(6.2,-1.5){\shortstack{$1$}}
\put(6,-1.5){\circle*{.1}}
\put(5,0){\circle*{.2}}
 \put(4,2){\circle*{.2}}
 \put(4,1.5){\circle*{.2}}
\put(6,3){\circle*{.2}}
\put(6,2.5){\circle*{.2}}
\put(4,2){\line(2,1){2}}
\put(4,1.5){\line(2,1){2}} 
\put(0,-2){\line(1,0){12}}

\end{picture}
\end{center}
\bigskip
\bigskip
\bigskip

\begin{lema}{\rm  Let $G=y^m+b_1(x)y^{m-1}+\ldots+a_m(x)$ be a reduced
polynomial of ${\KK}((x^{1\over q}))[y]$ and let
$J=J(F,G)=F_xG_y-F_yG_x$ be the Jacobian of $F$ and $G$.  Let
$V(G)=\lbrace (\beta_{l_i},l_i), l_0 > l_1
>\ldots> l_{v_G}\rbrace$ be the set of vertices of $N(G)$ and let
$E(G)=\lbrace{{\triangle}^G_1},\ldots,{{\triangle}}^G_{v_G}\rbrace$
be the set of edges of $N(G)$. Assume that the following holds:

i) $k_{v_F}=l_{v_G}=0$, $\alpha_{v_F}\not=0$ and
$\beta_{v_G}\not=0$, i.e. $N(F)$  and $N(G)$ meet  the $x$-axis into
points different from the origin.

ii) $(\alpha^1,1)\in {\rm Supp}(F_{\triangle^F_{v_F}})$ (resp.
$(\beta^1,1)\in{\rm Supp}(G_{{\triangle^G_{v_G}}})$)  for some
$\alpha^1$ (resp. $\beta^1$) in $\mathbb{Q}$, and
$(\alpha^1,1)\notin V(F)$ (resp. $(\beta^1,1)\notin V(G))$.

iii) max$(O(F))> {\rm max}(O(G))$.

\noindent Then we have:

i) max$(O(J))={\rm max}(O(F_y))={\rm max}(O(F))$.

ii) If  $G_{{\triangle^G_{v_G}}}(x,0)=ax^{\beta_{l_{v_G}}},
a\in\KK^*$, then
 $(\alpha^1+\beta_{l_{v_G}}-1,0)\in V(J)$ and $J_{\triangle^J_{v_J}}=(-F_yG_x)_{\triangle^{F_yG_x}_{v_{F_yG_x}}}=-a\beta_{l_r}.x^{\beta_{l_r}-1}(F_y)_{\triangle^{F_y}_{v_F}}$.

}
\end{lema}

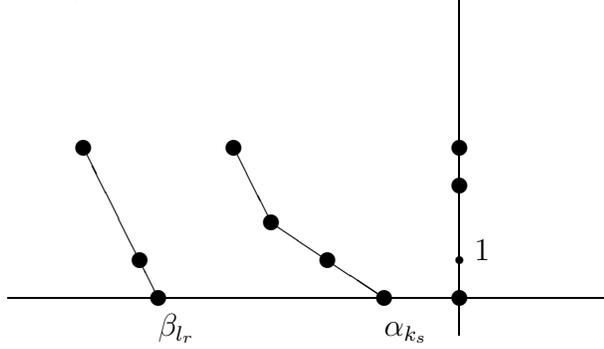
\begin{center}
  \unitlength=1cm
  \begin{picture}(6,4)(3,-2.5)
    \put(6,-2.5){\line(0,10){4.5}} 

    \put(3,0){\circle*{.2}}
    \put(6,0){\circle*{.2}}

    \put(6,-0.5){\circle*{.2}}

    \put(3,0){\line(1,-2){0.5}}
    \put(3.5,-1){\circle*{.2}}

    \put(3.5,-1){\line(3,-2){1.5}}
    \put(2,-2){\line(-1,2){1}}
    \put(2,-2.5){\shortstack{$\beta_{l_r}$}}
    \put(2,-2){\circle*{.2}}
    \put(1.75,-1.5){\circle*{.2}}
    \put(1,0){\circle*{.2}}
    \put(5,-2){\circle*{.2}}
    \put(5,-2.5){\shortstack{$\alpha_{k_s}$}}
    \put(4.25,-1.5){\circle*{.2}}

    \put(6,-2){\circle*{.2}}
    \put(6.2,-1.5){\shortstack{$1$}}
    \put(6,-1.5){\circle*{.1}}

    \put(0,-2){\line(1,0){8}}

  \end{picture}
\end{center}

\begin{demostracion}{.} It follows from the hypotheses that $(\alpha_{v(F)}-1,0)\in V(F_x)$, $(\alpha^1,0)\in V(F_y)$,  $(\beta_{v(G)}-1,0)\in V(G_x)$, and $(\beta^1,0)\in V(G_y)$. In particular
$(\alpha_{v(F)}+\beta^1-1,0)\in V(F_xG_y)$ and
$(\beta_{v(G)}+\alpha^1-1,0)\in V(F_yG_x)$. Since max$(O(F)) > {\rm
  max}(O(G))$, then
$\beta_{l_{v_G}}-\beta^1<\alpha_{k_{v(F)}}-\alpha^1$, in particular
$\beta_{l_{v_G}}+\alpha^1-1<\alpha_{k_{v_F}}+\beta^1-1$, and
$(\beta_{l_{v_G}}+\alpha^1-1,0)\in V(J)$. A similar argument shows
that the last edge of $J=F_xG_y-F_yG_x$ is nothing but the last edge
of $-F_yG_x$, and that
$(-F_yG_x)_{\triangle^{F_yG_x}_{v_{F_yG_x}}}=-a\beta_{l_r}.x^{\beta_{l_r}-1}(F_y)_{\triangle^{F_y}_{v_F}}$.
$\blacksquare$

\end{demostracion}

\section{Deformation of Newton polygons and applications}

\medskip

\noindent Let $f=y^n+a_{1}(x)y^{n-1}+\ldots+a_{n-1}(x)y+a_n(x)$ be a
reduced monic polynomial of $\KK((x))[y]$ and let ${\rm
Root}(f)=\lbrace y_1,\ldots,y_n\rbrace$. Let
 $f_1,\ldots,f_{\xi(f)}$ be the set of irreducible components of $f$ in
$\KK((x))[y]$.


\begin{definicion} {\rm Let $N$ be a nonnegative integer and let
$\gamma(x)=\sum_{k\geq k_0}a_kx^{k\over N}\in\KK((x^{1\over N}))$. Let
$M$ be a real number. We set




\[
 \gamma_{< M}=
\begin{cases}

\displaystyle{\sum_{k\geq k_0, {k\over N} < M}}a_kx^{k\over
N}&\text{if $M> \displaystyle{{k_0\over N}}$} \\
0&\text{otherwise}
\end{cases}
\]

\noindent and we call $\gamma_{< M}$ the $< M$-truncation of $\gamma(x)$.

\noindent Let $\theta$ be a generic element of ${\bf K}$. We
set

\[
 \gamma_{< M,\theta}=
\begin{cases}

\displaystyle{\sum_{k\geq k_0, {k\over N} < M}}a_kx^{k\over
N}+\theta.x^M &\text{if $M\geq \displaystyle{{k_0\over N}}$} \\
\theta.x^M&\text{otherwise}
\end{cases}
\]

\noindent and we call $\gamma_{< M, \theta}$ the $M$-deformation of
$\gamma(x)$.
 }\end{definicion}

\noindent Let $N$ be a nonnegative integer and let $\gamma(x)\in
\KK((x^{1\over N}))$. Let $M$ be a real number and let $\gamma_{<M}$
be the $<M$-truncation  of $\gamma(x)$. Consider the change of
variables $X=x, Y=y-\gamma_{< M}$. The polynomial
$F(X,Y)=f(X,Y+\gamma_{< M})$
is a monic polynomial of
degree $n$ in $Y$  whose coefficients are fractional meromorphic series
in $X$. Let $V(F)=\lbrace P_i=(\alpha_{k_i},k_i)|
i=1,\ldots,v_F\rbrace$ and let $E(F)=\lbrace
\triangle^F_1,\ldots,\triangle^F_{v_F}\rbrace$.

\begin{lema}{\rm  Let the notations be as above. Assume that   $\gamma\notin{\rm Root}(f)$ and let $M={\rm max}_{1\leq j\leq n}{\rm
c}(\gamma,y_j)$. We have the following:

i) ${\rm Root}(F(X,Y))=\lbrace
Y_k=y_k-\gamma_{< M}, k=1,\ldots,n\rbrace$.

ii) $O(F)=\lbrace {\rm c}(y_k,\gamma)| k=1,\ldots,n\rbrace.$

iii) There is exactly $k_i-k_{i+1}$ roots $y(x)$ of $F$ whose
contact with $\gamma$ is
$\dfrac{\alpha_{i}-\alpha_{i-1}}{k_i-k_{i-1}}$.

iii) The initial coefficients of Root$(F)$, denoted inco$(F)$, is $=\lbrace {\rm inco}(y_k-\gamma)| k=1,\ldots,n \rbrace.$

\noindent In particular, the Newton polygon $N(F)$ gives us a
complete information about the relationship between $\gamma(x)$ with
the roots of $f$. We call it the  Newton polygon of $f$ with respect
to $\gamma(x)$, and we denote it by $N(f,\gamma)$.}\end{lema}

\begin{demostracion}{.} We have

$$
F(X,Y)=f(X,Y+\gamma(X))=\prod_{k=1}^n(Y+\gamma(X)-y_k(X))
$$

\noindent now use Lemma 3.1.$\blacksquare$
\end{demostracion}

\begin{lema}{\rm Let $y_i(x)$ be a root of $f(x,y)=0$ and let

$$
M={\rm max}_{j\not= i}c(y_i,y_j).
$$

\noindent  Let $\tilde{y_i}={y_i}_{< M,\theta}=(y_i)_{<M}(x)+\theta x^M$ be the
$M$-deformation of $y_i$
 and consider the change of variables
$X=x,Y=y-\tilde{y_i}(X)$. Let $F(X,Y)=f(X,Y+\tilde{y_i}(X))$. We have the following:

i) $O(F)=\lbrace {\rm c}(y_j-y_i)| j\not=i\rbrace \rbrace.$

ii) $M={\rm max}(O(F))$.

iii) The last vertex of $N(F)$ belongs to the $x$-axis.

iv) Let $\triangle^F_{v_F}$ be the last edge of $N(F)$. We have
$(\alpha^1,1)\in{\rm Supp}(F_{\triangle^F_{v_F}})$ for some $\alpha^1$. Furthermore,
$(\alpha^1,1)\notin V(F)$.}
\end{lema}

\begin{demostracion}{.} We have

$$
F(X,Y)=f(X,Y+\tilde{y_i}(X))=\prod_{k=1}^n(Y+(y_i)_{<M}(X)+\theta X^M-y_k(X))
$$

\noindent and by hypothesis, $O((y_i)_{<M}(X)+\theta
X^M-y_k(X))=O(y_i(X)-y_k(X))$ for all $k\not= i$. Furthermore, $O((y_i)_{<M}(X)+\theta
X^M-y_i(X))=M=O(y_i(X)-y_j(X))$ for some $j\not= i$. This implies i)
and ii). Now $F(X,0)=\prod_{k=1}^n((y_i)_{<M}(X)+\theta
X^M-y_k(X))\not= 0$, hence iii) follows. Let $\triangle^F_{v_F}$ be
the last edge of $N(F)$ and let $y_{j_1},\ldots,y_{j_p}$ be the set of
roots of $f$ such that $c(y_i-y_{j_k})=M$ for all
$k=1,\ldots,l$. Write $y_i=\sum_{p}c^i_px^p$ and let
$y_i-y_{j_k}=c_{a_k}x^M+...$ for all $k=1,\ldots,l$. It follows that
$(Y_i)_{<M}(X)+\theta x^M-y_{j_k}(X)=(c_{a_k}+\theta)x^M+...$. Finally
$F_{\triangle^F_{v_F}}=(y-(c_M+\theta)x^M)\prod_{k=1}^l(y-(c_{a_k}+\theta)x^M)$. Since
$\theta$ is generic and $l\geq 1$, then iv) follows immediately.$\blacksquare$
\end{demostracion}

\noindent In particular, using the results of Section 3., the last
vertex of $N(F_Y)$ is $(\alpha^1,0)$,
$O(F)=O(F_y)$, and max$(O(F_Y))=M$. But
$F_Y(X,Y)=f_y(X,Y+\tilde{y_i}(X))$. This with the above Lemma led
to the following Proposition (see also [8], Lemma 3.3.):

\begin{lema} {\rm For $y_i(x),y_j(x), i\not=j$, there is a root $z_k(x)$
of $f_y(x,y)=0$ such that

$$
c(y_i(x),y_j(x))= c(y_i(x),z_k(x))=c(y_j(x),z_k(x)).
$$

\noindent Conversely, given $y_i(x),z_k(x)$, there is $y_j(x)$ for
which the above equality holds. Moreover, given $y_i(x)$ and $M
\in \mathbb{R}$,

$$
{\rm card} \lbrace y_j(x)|c(y_i(x),y_j(x))=M\rbrace= {\rm card}  \lbrace
z_k(x)| c(y_i(x),z_k(x))=M \rbrace.
$$
}\end{lema}

\begin{demostracion}{.} Let $i\not=j$ and let $M=c(y_i-y_j)$. Let
  $\tilde{y_i}=(y_i)_{<M}+\theta x^M$ be the $M$-deformation of $y_i$. Consider, as in
  Lemma 4.3., the change of variables $X=x,Y=y-\tilde{y_i}(X)$ and let
  $F(X,Y)=f(X,Y+\tilde{y_i}(X))$. It follows from Lemma 4.3. that
  $F(X,0)\not=0$, and if ${\rm deg}_yF_{\triangle^F_{v_F}}=r+1$, then
  there is $r$ roots $y_{j_1},\ldots,y_{j_r}$ of $f(x,y)=0$ such that
$c(y_i-y_{j_k})=M$ for all $k=1,\ldots,r$. Since $(\alpha^1,0)\in {\rm
  Supp}(F_{\triangle^F_{v_F}})$ for some $\alpha^1$, then the
cardinality of $E(F_y)$ is the same as the cardinality of
$E(F)$. Furthermore, $N(F_y)$ is
a translation of $N(F)$ with respect of the vector $(0,-1)$. Finally,
$(F_y)_{\triangle^F_{v_F}}=(F_{\triangle^F_{v_F}})_y$ is a
polynomial of degree $r$ in $y$. In particular, by Lemma 4.3., there
is $r$ roots of $f_y(x,y)=0$ whose contact with $y_i$ is $M$. This
completes the proof of the result.$\blacksquare$
\end{demostracion}

\noindent Let $g=y^m+a_1(x)y^{m-1}+\ldots+a_m(x)$ be a reduced
monic polynomial of ${\KK}((x))[y]$ and denote by
$z_1,\ldots,z_m$ the set of roots of $g$. Let
$y_i(x)\in{\rm Root}(f)$ and let:

$$
M={\rm max}(\lbrace c(y_i,y_j)|j\not= i\rbrace\cup \lbrace
c(y_i,z_k)|k=1,\ldots,m\rbrace)
$$

\begin{lema}{\rm Let the notations be as above, and assume  that $M >
    {\rm max}_{1\leq k\leq m}c(y_i,z_k)$. Let $\tilde{y_i}=(y_i)_{<M}+\theta x^M$ be
    the $M$-deformation
of $y_i$ and consider the change of variables $X=x,
Y=y-\tilde{y_i}(X)$. Let $F(X,Y)=f(X,Y+\tilde{y_i}(X))$,
$G(X,Y)=g(X,Y+\tilde{y_i}(X))$. We have the following

i) $F(X,0)\not=0$ and $G(X,0)\not=0$, i.e. $N(F)$ and $N(G)$
meet the $x$-axis.

ii) max$(O(F))=M > {\rm max}(O(G))$

iii) If $\triangle^F_{v_F}$ (resp. $\triangle^G_{v_G}$) denotes the
last edge of $N(F)$ (resp. $N(G)$) then $(\alpha^1,1)\in{\rm
Supp}(F_{\triangle_{v_F}})$ (resp. $(\beta^1,1)\in{\rm
Supp}(G_{\triangle_{v_G}})$)  for some $\alpha^1$
(resp. $\beta^1$), and $(\alpha^1,1)\notin V(F)$ (resp.
$(\beta^1,1)\notin V(G)$).}
\end{lema}

\begin{demostracion}{.} Let $F(X,Y)=\prod_{j=1}^n(Y-Y_j(X))$ and
  $G(X,Y)=\prod_{k=1}^m(Y-Z_k(X))$. Clearly
  $Y_j(X)=y_j(X)-(y_i)_{<M}(X)+\theta X^M$,
  $Z_k(X)=z_k(X)-(y_i)_{<M}(X)+\theta X^M$. In particular,  for all
  $k=1,\ldots,m$, $O(Z_k)=c(y_i,z_k)<M$. On the other hand, for all
  $j\not= i$, $O(Y_j)=c(y_i,y_j)\leq M$ with equality for at least one
  $j$, and $O(Y_i)=M$. This implies i) and ii). Now iii) follows by a
  similar argument as in Lemma 4.3.$\blacksquare$

\end{demostracion}

\noindent Let $J=J(f,g)$, and note that $J(F,G)=J(X,Y)$. In particular, by the results of the
previous Section we get the following:

\begin{lema} {\rm For $y_i(x),y_j(x), i\not=j$, if $c(y_i,y_j)> {\rm max}_{1\leq k\leq m}c(y_i,z_k)$,
then there is a root $u_l(x)$ of $J(x,y)=0$ such that

$$
c(y_i(x),y_j(x))= c(y_i(x),u_l(x))
$$

\noindent Conversely, given $y_i(x),u_l(x)$, if $c(y_i,u_l)>
c(y_i,z_k),k=1,\ldots,m$, there is $y_j(x)$ for which the above
equality holds. Moreover, given $y_i(x)$ and $M \in \mathbb R$, if
$M >{\rm max}_{1\leq k \leq m}c(y_i,z_k)$, then:

$$
{\rm card} \lbrace y_j(x)| c(y_i(x),y_j(x))=M\rbrace= {\rm card} \lbrace
u_l(x)| c(y_i(x),u_l(x))=M \rbrace.
$$
}\end{lema}

\begin{demostracion}{.} Let $M=c(y_i,y_j)$ and consider the change of
  variables $X=x, Y=y-\tilde{y_i}(X)$, where
  $\tilde{y_i}=(y_i)_{<M}+\theta x^M$ is the $M$-deformation of
  $y_i$. Let $F(X,Y)=f(X,Y+\tilde{Y_i}(X))$ and
  $G(X,Y)=g(X,Y+\tilde{Y_i}(X))$. Il follows from the hypotheses that
  $F$ and $G$ satisfies conditions i), ii), and iii) of Lemma 3.4. In
  particular
  $J(X,Y)_{\triangle^{J(X,Y)}_{v(J(X,Y)}}=(G_{\triangle^G_{v(G)}}(X,0))_X.(F_{\triangle^F_{v_F}})_Y$. The
      proof follows now by a similar argument as in Lemma 4.4.$\blacksquare$

\end{demostracion}










\section{Five main results}

\medskip

\noindent Let $f=y^n+a_1(x)y^{n-1}+\ldots+a_n(x)$ be a monic reduced
 polynomial of $\KK((x))[y]$ and let
$f=f_1.f_2.\ldots.f_{\xi(f)}$ be the decomposition of $f$ into
irreducible components of $\KK((x))[y]$. Let $f_y$ be the
$y$-derivative of $f$ and let ${\rm Root}(f)=\lbrace y_1(x),\ldots,y_n(x)\rbrace$.

\noindent For all $1\leq i\leq \xi(f)$,  set $n_{f_i}={\rm deg}_y(f_i)$,
and let $(m^{f_i}_k)_{1\leq k \leq h_{f_i}+1}, (d^{f_i}_k)_{1\leq
k\leq h_{f_i}+1}), (e^{f_i}_k)_{1\leq k \leq h_{f_i}}, \break
(r^{f_i}_k)_{1\leq k \leq h_{f_i}+1}$ be the set of characteristic
sequences of $f_i$. 

\begin {proposicion}{\rm  Assume that $\xi(f)=1$, i. e. $f=f_1$ is irreducible.  For all $1\leq k\leq
h_f$, we have:

 $$
 {\rm card} \lbrace z(x)\in{\rm Root}(f_y)|c(f,z(x))=\displaystyle{m^{f}_k\over n_{f}} \rbrace = (e^{f}_k-1)\displaystyle{n_{f}\over d^{f}_k}.
 $$}
 \end{proposicion}

\begin{demostracion}{.} Note that, by Lemma 4.4., $c(f,z(x))\in\lbrace
\displaystyle{{m^{f}_1\over n_{f}},\ldots, {m^{f}_{h_{f}}\over n_{f}}}\rbrace$.
Assume first that $k=h_{f}$ and fix a
root $y_p$ of $f$. By Lemma 1.2., $y_p$ has the contact
$\displaystyle{m^{f}_{h_{f}}\over n_{f}}$ with exactly
$d^{f}_{h_{f}}-d^{f}_{h_{f}+1}=d^{f}_{h_{f}}-1=e^{f}_{h_{f}}-1$
roots of $f$, consequently, by Lemma 4.4.,  there is
$e^{f}_{h_{f}}-1$ roots of $f_y$ whose contact with $y_p$ is
$\displaystyle{m^{f}_{h_{f}}\over n_{f}}$. Denote the set of these
roots by $D_p$. Each element of $D_p$ has the contact
$\displaystyle{m^{f}_{h_{f}}\over n_{f}}$ with exactly $d^{f}_{h_{f}}$
roots of $f$ (since we have to add $y_p$). Denote this set
by $C_p$. Let $y_q\not\in C_p$ be a root of $f$. Repeating  with $y_q$
what we did for $y_p$, we construct $D_q$ and $C_q$ in a
similar way. Obviously $C_p\cap C_q =\emptyset$ (otherwise,
$c(y_p,y_q)=\displaystyle{m^{f}_{h{f}}\over n_{f}}$, which is
impossible because $y_q\not\in C_p$). This implies that $D_p\cap D_q
=\emptyset$.... This process divides the $n_{f}$ roots of
$f$ into $\displaystyle{n_{f}\over d^{f}_{h_{f}}}$ disjoint groups
$C_1,\ldots,C_{{n_{f}\over d^{f}_{h_{f}}}}$ such that for all $1\leq
p\leq \displaystyle{n_{f}\over d^{f}_{h_{f}}}$, $C_p$ contains the roots
of $f$ having the contact $\displaystyle{m^{f}_{h_{f}}\over
n_{f}}$ with the elements of $D_p$. For all
$z(x)\in D_p, c(f,z(x))=\displaystyle{m^{f}_{h_{f}}\over
n_{f}}$, in particular

$$
  {\rm card} \lbrace z(x)\in{\rm Root}(f_y)| c(f,z(x))=\displaystyle{m^{f}_{h_{f}}\over n_{f}}\rbrace =
  \sum_{p
=1}^{n_{f}\over d^{f}_{h_{f}}} {\rm card}
D_p=(e^{f}_{h_{f}}-1)\displaystyle{n_{f}\over d^{f}_{h_{f}}}.
$$

\noindent Assume that the equality is true for $k=h_{f},\ldots,j+1$,
then there is exactly $\sum_{i =j+1}^{h_{f}}
(e^{f}_i-1)\displaystyle{{n_{f}\over d^{f}_j}=n_{f}-{n_{f}\over
d^{f}_j}}$ roots of $f_y$ having the  contact $\geq m^{f}_j$
with $f$. We now  repeat the same argument with
$\displaystyle{{n_{f}\over d^{f}_j}}$, $\displaystyle{d^{f}_j\over
d^{f}_{j+1}}-1$ instead of $n_{f}$ and
$d^{f}_{h_{f}}-1$.$\blacksquare$
\end{demostracion}

\begin{proposicion}{\rm  Let $M\in\mathbb{Q}$  and
let $1\leq i\leq \xi(f)$. Assume that  $M
\not=\displaystyle{m^{f_i}_{k}\over {n_{f_i}}}$ for all $k=1,\ldots, h_{f_i}$.
We have:

\begin{eqnarray*}
{\rm card} \lbrace z(x)\in{\rm Root}(f_y)|c(f_i,z(x))={M} \rbrace = {\rm card}  \lbrace y(x)\in{\rm Root}(f)|
c(f_i,y(x))=M\rbrace=
\sum_{c(f_i,f_k)={M}}n_{f_k}.
\end{eqnarray*}}
\end{proposicion}

\begin{demostracion}{.} Let, without loss
of generality, $i=1$  and let $k > 1$ be such  that
$c(f_1,f_k)= M$. Fix a root $y_p(x)$ of $f_1$. Since $c(y_p(x),f_k)=M$, then there is a root  $y(x)$ of $f_k$ such that
$c(y_p(x),y(x))= M$. Let $\theta\in\lbrace
0,\ldots,h_{f_1}\rbrace$ be the smallest integer such that $M <
\displaystyle{m^{f_1}_{\theta+1}\over {n_{f_1}}}$ and consider another root  $y_j(x)$ of
$f_1$. We have:

\[
c(y_j,y(x))=O_x(y_j-y(x))=O_x(y_j-y_p+y_p-y(x))=
\begin{cases}

{\displaystyle{M}}&\text{if $O_x(y_j-y_p)\geq {m^{f_1}_{\theta+1}\over n_{f_1}}$} \\
O_x(y_j-y_p) &\text{if $O_x(y_j-y_p)<{m^{f_1}_{\theta+1}\over
n_{f_1}}$}
\end{cases}
\]

\noindent By Lemma 1.2., there is exactly $d^{f_1}_{\theta+1}-1$
roots of $f_1$ having a contact
$\geq {\dfrac{m^{f_1}_{\theta+1}}{n_{f_1}}}$ with $y_p$, consequently, by the formula above,
there is exactly $d^{f_1}_{\theta+1}$ roots of $f_1$
having the contact $M$ with $y(x)$ (since
we have to add $y_p$).
Denote this set by $C_p$ and let $D_k^p$ be the set of roots
 of $f_k$ having the contact $M$ with
$y_p$. In particular an element  of $D_k^p$ has the contact
$M$ with every element of $C_p$.

\noindent Let $y_q\notin C_p$ be a root of $f_1$ and repeat
the same construction with $y_q$ instead of $y_p$. It is clear that
$C_p\cap C_q=\emptyset$ (otherwise, if $\bar{y}\in C_p\cap C_q$, then $c(\bar{y},y_p)=
c(\bar{y},y_q)=M$, in particular
 $c(y_p,y_q)\geq M$, which is a contradiction
 since $y_q\notin C_p$), in particular $D^p_k\cap
D^q_{k}=\emptyset$. This process divides the set of roots
of $f_k$ into disjoint $\displaystyle{n_{f_1}\over
d^{f_1}_{\theta+1}}$ groups $D_k^1,\ldots,D_k^{n_{f_1}\over
d^{f_1}_{\theta+1}}$: for all $1\leq p\leq \displaystyle{n_{f_1}\over
d^{f_1}_{\theta+1}}, D_k^{p}$ contains the roots of $f$ having
the contact $M$ with the elements of
$C_p$. Repeating what we did with another  $f_l$, $l\not= k$,
such that $c(f_1,f_l)=M$, then
adding the $D_k^{p}$'s, We obtain  disjoint
$\displaystyle{n_{f_1}\over d^{f_1}_{\theta+1}}$ groups
$D_1,\ldots,D_{n_{f_1}\over d^{f_1}_{\theta+1}}$ such that
$D_p$ contains the roots of $f$ having the contact
$M$ with the elements of $C_p$. We have,
by Lemma 4.4.

$$
{\rm card} \lbrace z(x)\in{\rm Root}(f_y)|c(y_p(x),z(x))=M\rbrace={\rm card} \lbrace y(x)\in{\rm Root}(f)|
c(y_p(x),y(x))=M \rbrace = {\rm card} D_p
$$

\noindent Let $z(x)\in{\rm Root}(f_y)$ and assume that
c$(z(x),y_p)=M$. If $y_q\in{\rm Root}(f), y_q\not= y_p$,
since $c(y_p,y_q)\not=M$,
then $c(z(x),y_q)\leq M$. In particular
c$(f_1,z(x))=M$. Finally

 $$
 {\rm card} \lbrace z(x)\in{\rm Root}(f_y)| c(f_1,z(x))=M
 \rbrace= \sum_{p=1}^{{n_{f_1}\over d^{f_1}_{\theta+1}}}{\rm card} D_p =\sum_{c(f_1,f_k)=M}n_{f_k}
 $$

\noindent This proves our assertion.$\blacksquare$
 \end{demostracion}

\medskip

\begin {proposicion}{\rm  For all  $1\leq i\leq r$ and for all $1\leq
\theta\leq h_{f_i}$, we have:

$$
{\rm card} \lbrace z(x)\in{\rm Root}(f_y)|c(f_i,z(x))=\displaystyle{m^{f_i}_\theta\over
n_{f_i}} \rbrace ={\rm card} \lbrace y(x)\in{\rm Root}(f)|c(f_i,y(x))=
\displaystyle{m^{f_i}_\theta\over n_{f_i}}\rbrace +
(e^{f_i}_\theta-1)\displaystyle{n_{f_i}\over d^{f_i}_\theta}
$$

$$
\quad\quad\quad\quad\quad\quad\quad
=\sum_{c(f_i,f_k)={m^{f_i}_\theta\over
n_{f_i}}}n_{f_k}+(e^{f_i}_\theta-1)\displaystyle{n_{f_i}\over
d^{f_i}_\theta}.
$$}

\end{proposicion}

\begin{demostracion}{.} Let, without loss of generality, $i=1$ and
assume that $\displaystyle{c(f_1,f_k)={m^{f_1}_\theta\over n_{f_1}}}$ for at least
one $k>1$. Let $y_p$ be a root of $f_1$. Since
$c(f_k,y_p)=\displaystyle{m^{f_1}_{\theta}\over n_{f_1}}$, then there
is  a root $y(x)$ of $f_k$ such that
$c(y_p(x),y(x))=\displaystyle{m^{f_1}_{\theta}\over n_{f_1}}$. By Lemma
1.2., there is exactly $d^{f_1}_{\theta}-1$ roots
of $f_1$ having a contact $\geq \displaystyle{m^{f_1}_{\theta}\over
n_{f_1}}$ with $y_p$. Let $y_j(x)$ be a root of $f_1$ such
that $c(y_p,y_j)\geq \displaystyle{m^{f_1}_{\theta}\over n_{f_1}}$,
then $c(y_j,y(x))= O_x(y_j-y(x))=O_x(y_j-y_p+y_p-y(x))\geq
\displaystyle{m^{f_1}_{\theta}\over {n_{f_1}}}$. On the other hand,
$c(y_j(x),y(x))\leq c(f_1,f_k)=\displaystyle{m^{f_1}_{\theta}\over
n_{f_1}}$, hence $c(y_j,y(x))=\displaystyle{m^{f_1}_{\theta}\over
n_{f_1}}$. Consequently, there is exactly $d^{f_1}_{\theta}$
roots of $f_1$ having the contact
$\displaystyle{m^{f_1}_{\theta}\over n_{f_1}}$ with $y(x)$. Denote
this set by $C_p$ and let $D^k_p$ be the set of roots of
$f_k$ such that for all $y(x)\in D^k_p,
c(y_p(x),y(x))=\displaystyle{m^{f_1}_{\theta}\over n_{f_1}}$. In
particular, an element  of $D^k_p$ has the contact
$\displaystyle{m^{f_1}_{\theta}\over n_{f_1}}$ with every element of
$C_p$.

\noindent Let $y_q\notin C_p$ be a root of $f_1$ and repeat
the same construction with $y_q$ instead of  $y_p$. We have,
by a similar argument as in Proposition 5.2.,  $C_p\cap
C_q =\emptyset$ and consequently $D^k_p\cap D^k_q=\emptyset$. This
divides the set of roots of $f_k$ into disjoint
$\displaystyle{n_{f_1}\over d^{f_1}_{\theta}}$ groups
$D_k^1,\ldots,D_k^{n_{f_1}\over d^{f_1}_{\theta}}$. Each element of
$D_k^{p}$ has the contact
$\displaystyle{m^{f_1}_{\theta}\over n_{f_1}}$ with the elements of
$C_p$. Repeating the same argument with the set of $f_l$ such that
$c(f_1,f_l)=\displaystyle{m^{f_1}_{\theta}\over n_{f_1}}$,
then adding the $D^k_p$'s, we obtain  disjoint
$\displaystyle{n_{f_1}\over d^{f_1}_{\theta}}$ groups
$D_1,\ldots,D_{n_{f_1}\over d_{\theta}}$ such that for all $1\leq
p\leq \displaystyle{n_{f_1}\over d_{\theta}}$, $D_p$ contains the
roots of $\displaystyle{f\over f_1}$ having the contact
$\displaystyle{m^{f_1}_{\theta}\over n_{f_1}}$ with the elements of
$C_p$. We have, by Lemma 2.2. and Lemma 4.4.:

$$
{\rm card} \lbrace z(x)\in{\rm Root}(f_y)| c(y_p,z(x))=\displaystyle{m^{f_1}_{\theta}\over
n_{f_1}} \rbrace ={\rm card} \lbrace
y(x)|c(y_p,y(x))=\displaystyle{m^{f_1}_{\theta}\over n_{f_1}}\rbrace
$$

$$
={\rm card} D_p + (e^{f_1}_{\theta}-1)\displaystyle{n_{f_1}\over
d^{f_1}_{\theta}}
$$

\noindent And  by a similar argument as in Proposition 5.2.,

 $$
 {\rm card} \lbrace z(x)\in{\rm Root}(f_y)|c(f_1,z(x))=\displaystyle{m^{f_1}_{\theta}\over n_{f_1}}\rbrace=
 (\sum_{1\leq p\leq {n_{f_1}\over d^{f_1}_{\theta}}}
 {\rm card} D_p)+(e^{f_1}_{\theta}-1)\displaystyle{n_{f_1}\over d^{f_1}_{\theta}}
 $$
 $$
 \quad\quad\quad\quad\quad\quad\quad\quad\quad\quad\quad\quad\quad\quad\quad\quad\quad=
 (\sum_{c(f_1,f_k)={m^{f_1}_{\theta}\over n_{f_1}}}n_{f_k})+(e^{f_1}_{\theta}
-1)\displaystyle{n_{f_1}\over d^{f_1}_{\theta}}
$$

\noindent This proves our assertion.$\blacksquare$
\end{demostracion}

\noindent Let $g=y^m+b_1(x)y^{m-1}+\ldots+b_m(x)$ be a  monic reduced
polynomial of
${\KK}((x))[y]$ and let $g_1,\ldots,g_{\xi(g)}$ be the set of irreducible components
of $g$ in ${\KK}((x))[y]$. Let Root$(g)=\lbrace z_1,\ldots,z_m\rbrace$, and
let $J=J(f,g)$ be the Jacobian of $f$ and
$g$.

 \begin{proposicion}{\rm  Let $M\in \mathbb{Q}$ and assume
 that $c(y(x),y'(x))=M$ for some $y(x),y'(x)\in{\rm Root}(f)$, and
  that $M >{\rm max}_{1\leq j\leq m}c(y(x),z_j(x))$. Let $1\leq i\leq \xi(f)$, and assume
  that  $M\not=\displaystyle{m^{f_i}_{k}\over{n_{f_i}}}$ for all $k=1,\ldots, h_{f_i}$.
We have the following

\begin{eqnarray*}
{\rm card} \lbrace u(x)\in{\rm
  Root}(J)|c(f_i,u(x))=M \rbrace = {\rm card}  \lbrace y_j(x)|
c(f_i,y_j(x))=\displaystyle{M}\rbrace=
\sum_{c(f_i,f_k)=M}n_{f_k}.
\end{eqnarray*}}
\end{proposicion}
\begin{demostracion}{.} The proof is similar to the proof of Proposition 5.2., where
Lemma 4.4. is replaced by Lemma 4.6.$\blacksquare$
\end{demostracion}

 \begin{proposicion}{\rm  Let $1\leq \theta\leq h_{f_i}$ and assume
 that $\displaystyle{m^{f_i}_\theta\over
n_{f_i}}> {\rm max}_{1\leq j\leq n, 1\leq k\leq m}c((y_j(x),z_k(x))$. We have the following

$$
{\rm card} \lbrace u(x)\in{\rm Root}(J)|c(f_i,u(x))=\displaystyle{m^{f_i}_\theta\over
n_{f_i}} \rbrace ={\rm card} \lbrace y(x)\in{\rm Root}(f)|c(f_i,y(x))=
\displaystyle{m^{f_i}_\theta\over n_{f_i}}\rbrace +
(e^{f_i}_\theta-1)\displaystyle{n_{f_i}\over d^{f_i}_\theta}
$$

$$
\quad\quad\quad\quad\quad\quad\quad
=\sum_{c(f_i,f_k)={m^{f_i}_\theta\over
n_{f_i}}}n_{f_k}+(e^{f_i}_\theta-1)\displaystyle{n_{f_i}\over
d^{f_i}_\theta}.
$$}
\end{proposicion}
\begin{demostracion}{.} The proof is similar to the proof of Proposition 5.3., where
Lemma 4.4. is replaced by Lemma 4.6.$\blacksquare$
\end{demostracion}

\section{The tree of contacts}

\medskip
\noindent Let
$f$ be a monic reduced polynomial  in
$\KK((x))[y]$ and let $f=f_1.\ldots.f_{\xi(f)}$ be the factorization of
$f$ into irreducible components of $\KK((x))[y]$. We define the set of contacts of $f$ to be the set:

$$
C(f)= \lbrace c(f_p,f_q)| 1\leq p\not=q\leq \xi(f)\rbrace\cup
\cup_{k=1}^{\xi(f)} \lbrace {m^{f_k}_1\over
n_{f_k}},\ldots,{m^{f_k}_{h_{f_k}}\over n_{f_k}}\rbrace
$$

 \noindent Let $C(f)=\lbrace M_1,\ldots,M_{t_f}\rbrace$. The tree
 associated with $f$ is constructed as follows:

\medskip

\noindent Let $M\in C(f)$ and define  $C_M(f)$ to be the set of irreducible
components  of $f$ such that

 $$
f_p\in C_M(f)\Leftrightarrow c(f_p,f_q)\geq
M \text{ for some } 1\leq q\leq \xi(f)
$$

\noindent with the understanding that $c(f_k,f_k)\geq M$ if and only if
$\displaystyle{m^{f_k}_i\over {n_{f_k}}}\geq M$ for some $1\leq i\leq h_{f_k}$. We associate
with $M$ the equivalence
relation on the set $C_M(f)$, denoted
$R_M$, and defined as follows:

$$
f_pR_M f_q\text{ if and only if } c(f_p,f_q) \geq M.
$$

\noindent We define the points of the tree $T(f)$ at the level $M$
to be the set of equivalence classes of $R_M$, and we denote this
set by $P_1^M,\ldots,P_{s_M}^M$. We shall say that $P_j^N$ dominates
$P_i^M$, and we write $P_j^N\geq P_i^M$,  if $P_j^N\subseteq P_i^M$.
 We shall say that $P_j^N$ strictly dominates
$P_i^M$, and we write $P_j^N > P_i^M$, if $P_j^N$ dominates $P_i^M$, $P_i^M\not=P_j^N$,
and  $C(f)\cap ]M,N[=\emptyset$. This defines an order on the set of
points of $T(f)$ with a unique minimal element, denoted  $P_1^{M_1}$. A point
$P_i^M$ is called a top point of $T(f)$ if it is maximal with
respect to this order. We denote by Top$(f)$ the set of top points of
$T(f)$.

\noindent Let $P_i^M$ be a point of $T(f)$, and let $P_{i_1}^{N_1},\ldots,P_{i_t}^{N_t}$
be the set of points that strictly dominate $P_i^M$. We set
$D_i^M=P_i^M-\cup_{l=1}^tP_{i_l}^{N_l}$. Clearly
$\lbrace P_{i_1}^{N_1},\ldots,P_{i_t}^{N_t},D_i^M\rbrace$ is a partition of $P_i^M$.
Furthermore, for all $F\in D_i^M$ and for all $F\not=G\in P_i^M, c(F,G)=M$. We also have the following:

i)  if $F\in D_i^M$, then $M\geq \displaystyle{m^F_{h_F}\over n_F}$

ii)   $P_i^M\in {\rm Top}(f)$ if and only if $P_i^M=D_i^M$.

\noindent If $P_j^N$ strictly dominates $P_i^M$, then
we link these two points be a segment of line. We define the set of
edges of  $T(f)$ to be the set of these  segments. Given a point $P_i^M$, if $D_i^M\not=\emptyset$, then
we associate with each $F\in D_i^M$ an arrow starting at the point $P_i^M$. Let
$P_1^{M_1},P_{i_2}^{M_2}\ldots,P_{i_k}^{M_k}$ be a set of points of $T(f)$
such that $P_{i_2}^{M_2}$ strictly dominates $P_1^{M_1}$,
$P_{i_k}^{M_k}\in{\rm Top}(f)$, and  $P_{i_j}^{M_j}$ strictly dominates
$P_{i_{j-1}}^{M_{j-1}}$ for all $3\leq j\leq k$. The union of edges
linking these points is called a branch of $T(f)$. Clearly, there are as many branches of $T(f)$
as there are top points of $T(f)$.

\vskip1cm
 \unitlength=1cm
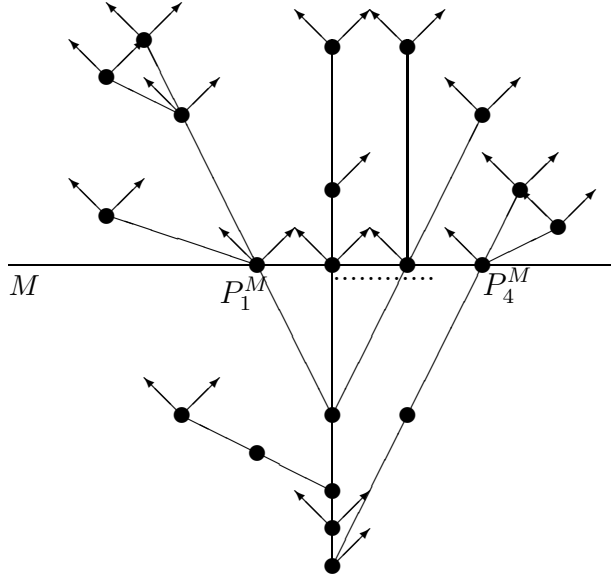
 \begin{center}
\begin{picture}(6,4)(3,-1)

\put(6,-4){\line(0,10){7}} 
\put(6,-4){\circle*{.2}} \put(6,-4){\vector(1,1){0.5}}
\put(6,-3.5){\circle*{.2}} \put(6,-3.5){\vector(1,1){0.5}}\put(6,-3.5){\vector(-1,1){0.5}}
\put(6,0){\circle*{.2}} \put(6,0){\vector(1,1){0.5}}\put(6,0){\vector(-1,1){0.5}}
\put(6,1){\circle*{.2}}\put(6,1){\vector(1,1){0.5}}
\put(6,2.9){\circle*{.2}}\put(6,2.9){\vector(1,1){0.5}}\put(6,2.9){\vector(-1,1){0.5}}

\put(6,-4){\line(1,2){2.5}}  
\put(7,-2){\circle*{.2}}
\put(8.5,1){\circle*{.2}}\put(8.5,1){\vector(1,1){0.5}}\put(8.5,1){\vector(-1,1){0.5}}

\put(6,-3){\line(-2,1){2}} 
\put(6,-3){\circle*{.2}} \put(5,-2.5){\circle*{.2}}
\put(4,-2){\circle*{.2}}\put(4,-2){\vector(1,1){0.5}}\put(4,-2){\vector(-1,1){0.5}}

\put(6,-2){\line(1,2){2}}  
\put(8,2){\circle*{.2}}
\put(8,2){\vector(1,1){0.5}}\put(8,2){\vector(-1,1){0.5}}

\put(6,-2){\circle*{.2}}
\put(6,-2){\line(-1,2){2.5}} 
\put(5,0){\circle*{.2}}
\put(5,0){\line(-3,1){2}}
\put(5,0){\vector(-1,1){0.5}}
\put(5,0){\vector(1,1){0.5}}
\put(3,0.65){\circle*{.2}}\put(3,0.65){\vector(1,1){0.5}}\put(3,0.65){\vector(-1,1){0.5}}

\put(4,2){\circle*{.2}}\put(4,2){\vector(1,1){0.5}}\put(4,2){\vector(-1,1){0.5}}

\put(3.5,3){\circle*{.2}} \put(3.5,3){\vector(1,1){0.5}}\put(3.5,3){\vector(-1,1){0.5}}

 \put(4,2){\line(-2,1){1}} 
\put(3,2.5){\circle*{.2}}  \put(3,2.5){\vector(1,1){0.5}}\put(3,2.5){\vector(-1,1){0.5}}

\put(7,0){\line(0,3){3}} 
\put(7,0){\circle*{.2}} \put(7,0){\vector(-1,1){0.5}}
\put(7,2.9){\circle*{.2}}\put(7,2.9){\vector(1,1){0.5}}\put(7,2.9){\vector(-1,1){0.5}}

\put(8,0){\line(2,1){1}}
\put(8,0){\circle*{.2}} \put(8,0){\vector(-1,1){0.5}}
\put(9,.5){\circle*{.2}} \put(9,.5){\vector(1,1){0.5}}\put(9,.5){\vector(-1,1){0.5}}

\put(1.7,0){\line(1,0){8}}
\put(1.7,-0.4){\shortstack{\small{$M$}}}
 \put(4.5,-0.5){\shortstack{$\small
P_1^M$}} \put(6,-0.2){\shortstack{$............$}}
\put(8,-0.4){$\displaystyle{\small{P_4^M}}$}
\end{picture}
\end{center}
\bigskip
\vskip3cm

\begin{lema}{\rm  Let $P_i^M$ be a  point of $T(f)$. We have the following:

i)  For all $F,G\in P_i^M, c(F,G)\geq M$.

ii) For all $F\in P_i^M$ and for all $H\notin P_i^M, c(F,H) < M$.

iii)  For all $F,G\in P_i^M$ and for all $H\notin P_i^M, c(F,H)=c(G,H)$. We denote
this rational by $c(H,P_i^M)$.

iv) let $F\in P_i^M$ and let $1\leq \theta\leq h_F+1$ be the smallest
integer such that $M \leq \displaystyle{m^F_{\theta}\over n_F}$.
If $\theta \geq 2$  then $\displaystyle{m^G_k\over
n_G}=\displaystyle{m^F_k\over n_F}$ for all $G\in P_i^M$ and for all $1\leq k\leq
\theta-1$. We denote this
rational number by $\displaystyle{m_k\over n}(P_i^M))$. As a
consequence $\displaystyle{n_G\over d^G_{k}}$ does not depend on $G\in P_i^M)$
 and $1\leq k\leq \theta$. We  denote this
 rational number
 by $\displaystyle{n\over d_k}(P_i^M)$.

}
\end{lema}
\begin{demostracion}{.} The proof is an easy application of Lemma 1.5. and Lemma 1.6.$\blacksquare$
\end{demostracion}

\noindent Let  $P_i^M$ be a point of $T(f)$ and define the
subsets $X_1(M,i),\ldots,X_{s(M,i)}(M,i)$ of $P_i^M$ as follows:

- For all $k$ and for all $F,G \in X_k(M,i)$, $c(F,G)=M$.

- Given $F\in X_k(M,i)$ and $l\not=k$, if $F\notin X_l(M,i)$, then
$c(F,G)> M$ for some $G\in X_l(M,i)$ (in particular
$F,G\in P_j^N$ for some $P_j^N > P_i^M$).

\noindent The sets defined above satisfy the following property:

\begin{lema}{\rm  The cardinality of $X_k(M,i)$ does not depend on
  $1\leq k\leq s(M,i)$. We denote this cardinality by $c(M,i)$.}\end{lema}

\begin{demostracion}{.} Assume that $s(M,i)\geq 2$ and
let $1\leq a\not= b\leq s(M,i)$. We shall construct a bijective map
from $X_a(M,i)$ to $X_b(M,i)$. Let $F\in X_a(M,i)$.  If $F\notin
X_b(M,i)$, then there is $\tilde{F}\in X_b(M,i)$ such that
$c(F,\tilde{F}) > M$. We claim that $\tilde{F}$ is the only element with this
property. In fact, if there is $\tilde{F}\not=G\in X_b(M,i)$ such that
$c(F,G))
> M$, then $M=c(\tilde{F},G)) \geq {\rm min}(c(F,\tilde{F}),c(F,G))> M$, which is
a contradiction. This defines a map

$$
\phi_{a,b}: X_a(M,i)\longmapsto X_b(M,i)
$$
\[
\phi_{a,b}(F)=
\begin{cases}

F&\text{if $F\in X_b(M,i)$} \\
\tilde{F}&\text{if $F\notin X_b(M,i)$}
\end{cases}
\]

\noindent This map is clearly bijective. This completes the proof of
the Lemma.$\blacksquare$\end{demostracion}

\begin{lema}{\rm Let the notations be as above, and let
$P_{i_1}^{N_{1}},\ldots,P_{i_t}^{N_{t}}$ be
the set of points that strictly dominate $P_i^M$. We have the following:

i) $D_i^M\subseteq X_k(M,i)$ for all $1\leq k\leq s(M,i)$.

ii) Given $1\leq k\leq s(M,i)$, $(X_k(M,i)\cap
P_{i_1}^{N_{1}},\ldots,X_k(M,i)\cap P_{i_t}^{N_{t}},
D_i^M)$ is a partition of $X_k(M,i)$.

iii) Given $1\leq k \leq s(M,i)$ and $1\leq l\leq t$,
$X_k(M,i)\cap P_{i_l}^{N_{1}}$ is reduced to one element.

iv) $c(M,i)=t+{\rm card}(D_i^M)$.}
\end{lema}

\begin{demostracion}{.} The first two assertions are clear,
on the other hand 3. $\Longrightarrow$ 4. We shall consequently
prove 3. Assume, without loss of generality, that $k=1$, and let $1\leq l\leq t$ .
Suppose that $X_1(M,i)\cap P_{i_l}^{N_{l}}= \emptyset$ and
let
$G\in P_{i_l}^{N_{l}}$. We have $c(F,G)=M$ for all  $F\in X_1(M,i)$,
in particular $G\in X_1(M,i)$, which is a contradiction. Consequently
$X_1(M,i)\cap P_{i_l}^{N_{l}}
\not= \emptyset$. Let  $G_1,G_2$ be two  polynomials of $X_1(M,i)\cap P_{i_l}^{N_{l}}$.
We have $c(G_1,G_2)=M$ and $c(G_1,G_2)\geq N>M$. This is a contradiction if
$G_1\not=G_2$.$\blacksquare$
\end{demostracion}

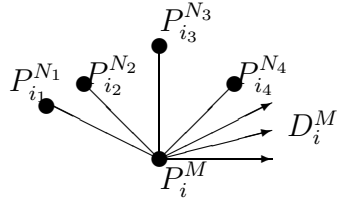
 \unitlength=1cm
 \begin{center}
\begin{picture}(6,4)(3,-1)

\put(6,0){\line(0,1){1.5}} 
\put(6,-0.4){$\displaystyle{P_{i}^{M}}$}
\put(6,0){\circle*{.2}}

\put(6,0){\line(1,1){1}}
\put(6,0){\line(-1,1){1}}
\put(6,0){\line(-2,1){1.5}}
\put(6,0){\vector(1,0){1.5}}
\put(6,0){\vector(4,1){1.5}}
\put(6,0){\vector(2,1){1.5}}

\put(7.7,0.3){$ D_i^M$}
\put(4,0.9){$\displaystyle{P_{i_1}^{N_1}}$}
\put(6,1.7){$\displaystyle{P_{i_3}^{N_3}}$}
 \put(7,1){\shortstack{$\small P_{i_4}^{N_4}$}}
\put(5,1){\circle*{.2}}
\put(4.5,0.7){\circle*{.2}}
\put(6,1.5){\circle*{.2}}
\put(7,1){\circle*{.2}}
\put(5,1){$\displaystyle{P_{i_2}^{N_2}}$}
\end{picture}
\end{center}

\noindent Let $P_i^M$ be a point of $T(f)$ and assume that $D_i^M\not=\emptyset$. For all $F\in D_i^M, c(F,F)\leq M$,
in particular  $M\geq {{m^F_{h_F}}\over {{n_F}}}$.

\begin{lema}{\rm Let the notations be as above. We have the following

\noindent i) If $M >{{m^{F}_{h_F}}\over {{n_F}}}$ for  all $F\in D_i^M$, then $n_F$
does not depend on $F\in D_i^M$. We denote it by $n(D_i^M)$. We have
$\sum_{F\in D_i^M}n_F=(c(M,i)-t)n(D_i^M)$.

\noindent ii) If $M= {{m^{F}_{h_F}}\over {{n_F}}}$ for some $F\in D_i^M$, then one of
the following hold

1$_{ii})$ $M={{m^{F}_{h_F}}\over {{n_F}}}$ for all $F\in D_i^M$. In this case, $n_F$
does not depend on $F\in D_i^M$. We denote it by $n(D_i^M)$. We have
$\sum_{F\in D_i^M}n_F=(c(M,i)-t)n(D_i^M)$.

1$_{iii)}$  $M>{{m^{F'}_{h_{F'}}}\over {{n_{F'}}}}$ for some $F'\in D_i^M$. In
this case, $M={{m^{F}_{h_F}}\over {{n_F}}}>{{m^{F'}_{h_{F'}}}\over {{n_{F'}}}}$ for all
$F\in D_i^M, F\not= F'$, and $n_F, h_F, (d^F_k)_{1\leq k\leq h_F}$ do not depend on
$F\in D_i^M, F\not= F'$. We denote these integers by $n(D_i^M),
h(D_i^M), d^{D_i^M}_k$. With these notations we have  $n(F')={\dfrac{n(D_i^M)}{d^{D_i^M}_{h(D_i^M)}}}$, and
$\sum_{F\in D_i^M}n(F)=(c(M,i)-t-1).n(D_i^M)+\dfrac{n(D_i^M)}{d^{D_i^M}_{h(D_i^M)}}$.
}
\end{lema}

\begin{demostracion}{.} By definition, for all $F,G\in D_i^M, c(F,G)=M$. Consequently
our results follow from Proposition 3.4.$\blacksquare$
\end{demostracion}




\noindent Let $H$ be a monic polynomial of ${\bf K}((x))[y]$ and
let $H_1,\ldots,H_{\xi(H)}$ be the set of irreducible components of
$H$ in $\KK((x))[y]$. Let $P_i^M$ be a point of $T(f)$ and let $F\in P_i^M$. We set:

$$
R_{=M}(F,H)=\prod_{c(F,H_k)=M}H_k
$$

\noindent and

$$
R_{>M}(F,H)=\prod_{c(F,H_k)>M}H_k
$$

\noindent In other words, $R_{=M}(F,H)$ (resp. $R_{>M}(F,H)$) is
the product of irreducible components of $H$ whose contact with
$F$ is $M$ (resp. $>M$).

\begin{lema}{\rm Suppose that $P_i^M\notin{\rm Top}(f)$ and
let $P_{i_1}^{N_1},\ldots,P_{i_t}^{N_t}$ be the set of points that strictly dominate
$P_i^M$. Fix $1\leq l\leq t$ and let $F\in P_{i_l}^{N_l}$. We have the following

i) For all
$G\in P_{i_l}^{N_l}$, $R_{=M}(G,H)=R_{=M}(F,H)$ (resp.
$R_{>M}(G,H)=R_{>M}(F,H)$). We denote this polynomial by
$R_{=M}(P_{i_l}^{N_l},H)$ (resp. $R_{>M}(P_{i_l}^{N_l},H)$).

\medskip

ii) For all $G\in P_{i_k}^{N_k}, k\not= l$, $R_{>M}(G,H)$ divides $R_{=M}(F,H)$.

\medskip

iii) For all $G\in D_i^M$, $R_{> M}(G,H)$ divides $R_{=M}(F,H)$.
}
\end{lema}

\begin{demostracion}{.} Let $\bar{H}$ be an irreducible component of $H$. If $G\in P_{i_l}^{N_l}$,
 then $c(F,G)\geq N_l > M$.
In particular, by Proposition 1.5., $c(G,\bar{H})=M$ (resp. $c(G,\bar{H}) > M$) if and only if
$c(F,\bar{H})=M$ (resp. $c(F,\bar{H})>M$). This proves i). If either $G\in P_{i_k}^{N_k}, k\not= l$ or
$G\in D_i^M$, then $c(F,G)=M$. In
particular, if $c(G,\bar{H})>M$, then, by Lemma 1.5.  $c(F,\bar{H})=M$. This proves ii) and iii).$\blacksquare$

\end{demostracion}

\noindent Let the notations be as above. It follows from ii), iii)
 of Lemma 6.5.  that:

$$
\prod_{k=2}^tR_{>M}(P_{i_k}^{N_k},H).\prod_{F\in
D_i^M}R_{>M}(F,H)\hskip0.08cm {\rm divides} \hskip0.08cm R_{=M}(P_{i_1}^{N_1},H).
$$

\noindent Set

$$
\overline{Q}_H(M,i)=\dfrac{R_{=M}(P_{i_1}^{N_1},H)}{\prod_{k=
2}^tR_{>M}(P_{i_k}^{N_k},H).\prod_{G\in D_i^M}R_{>M}(G,H)}
$$

\noindent and let

$$
Q_H(M,i)=\prod_{c(G,H_k)=M \forall G\in P_i^M}H_k
$$

\noindent i.e. $Q_H(M,i)$ is the product of the irreducible
components of $H$ whose contact with all $G\in P_i^M$ is $M$.

\begin{lema}{\rm With the notations above, we have
    $Q_H(M,i)=\overline{Q}_H(M,i)$.
}
\end{lema}

\begin{demostracion}{.} Let $\bar{H}$ be an irreducible component of $Q_H(M,i)$. For
all $G\in P_i^M, c(G,\bar{H})=M$. In particular, since $\cup_{k=1}^tP_{i_k}^{N_k}\subseteq P_i^M$, then
 $\bar{H}$ divides $R_{=M}(P_{i_1}^{N_1},H)$ and $\bar{H}$ does not divide ${\prod_{k=
2}^tR_{>M}(P_{i_k}^{N_k},H).\prod_{G\in D_i^M}R_{>M}(G,H)}$. Hence $Q_H(M,i)$
divides $\overline{Q}_H(M,i)$.

\noindent Let us prove that $\overline{Q}_H(M,i)$
divides ${Q}_H(M,i)$. Let $G\in P_i^M$ and
let   $\bar{H}$ be an irreducible component of
$\overline{Q}_H(M,i)$.

- If $G\in P_{i_1}^{N_1}$, then by Lemma 6.5. i), $R_{=M}(G,H)=R_{=M}(P_{i_1}^{N_1},H)$, in
particular $c(G,\bar{H})=M$.

- If $G\in  P_i^M- P_{i_1}^{N_1}$ then, by Lemma 6.3., $G\in D_i^M\cup
(\cup_{k=2}^tP_{i_k}^{N_k})$.  Suppose that $G\in D_i^M$. If $c(G,\bar{H}) > M$, then
$\bar{H}$ divides $R_{> M}(G,H)=R_{>M}(D_i^M,H)$. This contradicts the definition of
$\overline{Q}_H(M,i)$. In particular $c(G,\bar{H})=M$. By a similar argument we prove that
if $G\in \cup_{k=2}^tP_{i_k}^{N_k}$, then $c(G,\bar{H})=M$. This implies our assertion.$\blacksquare$
\end{demostracion}

\begin{lema}{\rm  Suppose that  $P_i^M\in {\rm Top}(f)$, and recall
    that in this case $P_i^M=D_i^M$. Let $F$ be an element of
    $D_i^M$. We have

$$
Q_H(M,i)=\dfrac{R_{=M}(F,H)}{\prod_{G\in D_i^M, G\not=
F}R_{>M}(G,H)}
$$
}
\end{lema}

\begin{demostracion}{.} The proof is similar to the proof of Lemma 6.6. $\blacksquare$
\end{demostracion}

\section{Factorization of the $y$-derivative}

\subsection{The irreducible case}
\medskip

\noindent  Let $f$ be a monic irreducible polynomial of
$\KK((x))[y]$ of degree $n_f$ in $y$ and consider the
characteristic sequences associated with $f$ as in Section 1. We
have the following:

\begin{proposicion} {\rm $f_y=P_1.\ldots.P_{h_f}$ and for all $k=1,\ldots,h_f$:

i) {\rm deg}$_yP_k=(e^f_k-1)\displaystyle{n_f\over d^f_k}$.

ii) {\rm int}$(f,P_k)=(e^f_k-1)r^f_k$.

iii) For all irreducible component $P$ of $P_k$,
$c(f,P)=\displaystyle{m^f_k\over n_f}$.}\end{proposicion}

\begin{demostracion}{.} i) and iii) result from Proposition 5.1.
and ii) results from Proposition 1.4.$\blacksquare$\end{demostracion}

\subsection {The general case}

\medskip

\noindent Let the notations be as in Section 5. In particular
$f$ is a monic reduced polynomial of $\KK((x))[y]$ and
$f_1,\ldots,f_{\xi(f)}$ are the irreducible components of $f$ in
$\KK((x))[y]$. Consider the characteristic sequences associated
with $f_1,\ldots,f_{\xi(f)}$ and let $T(f)$ be the tree of $f$. Fix a
point $P_i^M$ of $T(f)$.

\begin{lema}{\rm Let the notations be as above and let
    $P_i^M\in T(f)-{\rm Top}(f)$. If $D_i^M\not=\emptyset$, then deg$_y(R_{>M}(F,f_y))=0$ for all $F\in D_i^M$.
}
\end{lema}

\begin{demostracion}{.} Suppose that $P_i^M\notin {\rm Top}(f)$, and
  that $D_i^M\not= \emptyset$. Let $F\in D_i^M$. If
  deg$_y(R_{>M}(F,f_y))\not = 1$, then $c(F,H)=N>M$ for some irreducible
  component $H$ of $f_y$. In particular, by Lemma 4.4., $c(F,\bar{F})=N$
  for some irreducible component $\bar{F}$ of $f$, hence $F\in P_j^N$
  for some point $P_j^N\in T(f), N >M$. This is a contradiction
  because $F\in D_i^M$.$\blacksquare$
\end{demostracion}

\begin{lema}{\rm Suppose that
$P_i^M\notin {\rm Top}(f)$ and let
$P_{i_1}^{N_1},\ldots,P_{i_t}^{N_t}$ be the set of points of $T(f)$ that
strictly dominate $P_i^M$. We have:

$$
Q_{f_y}(M,i)=\dfrac{R_{=M}(P_{i_1}^{N_1},f_y)}{\prod_{l=2}^t R_{>
M}(P_{i_l}^{N_l},f_y)}
$$
}
\end{lema}

\begin{demostracion}{.} We have, by Lemma 6.6.:

$$
Q_{f_y}(M,i)=\dfrac{R_{=M}(P_{i_1}^{N_1},f_y)}{\prod_{l=2}^t R_{ >
M}(P_{i_l}^{N_l},f_y).\prod_{G\in D_i^M}R_{> M}(G,f_y)}
$$

\noindent On the other hand, by Lemma 7.2.,  if
$G\in D_i^M$, then deg$_y(R_{> M}(G,f_y)=0)$. This proves our
assertion.$\blacksquare$
\end{demostracion}

\noindent  Fix a polynomial $F_l\in
P_{i_l}^{N_l}$
for all $1\leq l\leq t$. By Lemma 6.6.,
$R_{=M}(P_{i_l}^{N_l},f_y)=R_{=M}(F_l,f_y)$
(resp. $R_{>M}(P_{i_l}^{N_l},f_y)=R_{>M}(F_l,f_y)$). In particular
we have:

$$
Q_{f_y}(M,i)=\dfrac{R_{=M}(F_1,f_y)}{\prod_{l=2}^t R_{>
M}(F_l,f_y)}
$$

\noindent The following Lemmas give the degrees of the two polynomials
describing $Q_{f_y}(M,i)$.

\begin{lema}{\rm Let $P_i^M$ be a point of $T(f)$ and let $\theta$ be the smallest integer such that $M\leq
\dfrac{m^F_{\theta}}{n_F}$ for all $F\in P_i^M$.  Let $(P_{i_l}^{N_l})_{1\leq l\leq t}$ be the set of points that strictly 
dominate $P_i^M$. Let $F_1\in P_{i_1}^{N_1}$. We have:

\small{
\[
{\rm deg}_yR_{=M}(F_1,f_y)=
\begin{cases}

\displaystyle{\sum_{l=2}^t(\sum_{F\in P_{i_l}^{N_l}}n_F)+\sum_{F\in D_i^M}n_F}&\text{if $M\not= {m_{\theta}^{F_1}\over {n_{F_1}}}$} \\
\displaystyle{\sum_{l=2}^t(\sum_{F\in P_{i_l}^{N_l}}n_F)+\sum_{F\in
D_i^M}n_F+(e^{F_1}_{\theta}-1){n_{F_1}\over
d^{F_1}_{\theta}}}&\text{if $M={m_{\theta}^{F_1}\over {n_{F_1}}}$}
\end{cases}
\]
}
}
\end{lema}

\begin{demostracion}{.} This results from  Propositions
5.2. and 5.3.$\blacksquare$
\end{demostracion}

\begin{lema}{\rm Let the notations and  the hypotheses
    by as in Lemma 7.4. We have:

$$
{\rm deg}_yR_{>M}(F_1,f_y)=\sum_{F\in P_{i_1}^{N_l}-\lbrace
F_1\rbrace}n_F+\sum_{M<{{m^{F_1}_j}\over n_{F_1}}}(e^{F_1}_j-1){n_{F_1}\over
d^{F_1}_j}
$$
}
\end{lema}
\begin{demostracion}{.} This results from  Propositions
5.2. and 5.3.$\blacksquare$
\end{demostracion}

\noindent As a corollary we have the following:

\begin{proposicion}{\rm Let the notations and the hypotheses
    by as in Lemma 7.4. and fix $F_l\in P_{i_l}^{N_{i_l}}$ for all $2\leq l\leq t$. We have:

\small{
\[
{\rm deg}_yQ_{f_y}(M,i)=
\begin{cases}

\displaystyle{\sum_{F\in
D_i^M}n_F+\sum_{l=2}^t[n_{F_l}}-\displaystyle{\sum_{M<{m^{F_l}_{j}\over {n_{F_l}}}}(e^{F_l}_j-1){n_{F_l}\over
d^{F_l}_j})]}&\text{if $M\not= {m^{F_1}_{\theta}\over n_{F_1}}$} \\
\displaystyle{\sum_{F\in D_i^M}n_F}+\sum_{l=2}^t[n_{F_l}-\sum_{M<  {m^{F_l}_j\over {n_{F_l}}}     }(e^{F_l}_j-1){n_{F_l}\over
d^{F_l}_j}]+(e^{F_1}_{\theta}-1){n_{F_1}\over d^{F_1}_{\theta}}&\text{if $M={m^{F_1}_{\theta}\over {n_{F_1}}}$}
\end{cases}
\]
}
}
\end{proposicion}

\begin{demostracion}{.} This results from Lemmas 7.4. and 7.5., since
  gcd($R_{=M}(F_1,f_y),R_{>M}(F_l,f_y))=1$ for all $2\leq l\leq t$.$\blacksquare$
\end{demostracion}

\noindent Note that with the hypotheses of Proposition 7.6.,
\small{
\[
\displaystyle{n_{F_l}-\sum_{M<   {m^{F_l}_{j} \over {n_{F_l}}}   }(e^{F_l}_j-1){n_{F_l}\over d^{F_l}_j}}=
\begin{cases}

\displaystyle{n_{F_l}-\sum_{j=\theta}^{h_{F_l}}(e^{F_l}_j-1){n_{F_l}\over
d^{F_l}_j}={n_{F_l}\over d^{F_l}_{\theta}}}&\text{if $M\not= {m^{F_l}_{\theta}\over n_{F_l}}$} \\
\displaystyle{n_{F_l}-\sum_{j=\theta+1}^{h_{F_l}}(e^{F_l}_j-1){n_{F_l}\over
d^{F_l}_j})={n_{F_l}\over d^{F_l}_{\theta+1}}}&\text{if
$M={m^{F_l}_{\theta}\over n_{F_l}}$}
\end{cases}
\]
}

\noindent  Let $A$ (resp. $B$) be the set of $1\leq l\leq t$ for
which $M = \displaystyle{m^{F_l}_{\theta}\over n_{F_l}} $ (resp. $M
<\displaystyle{m^{F_l}_{\theta}\over n_{F_l}}$). It follows that:
\small{
\[
{\rm deg}_yQ_{f_y}(M,i)=
\begin{cases}

\displaystyle{\sum_{l\in A}{n_{F_l}\over d^{F_l}_{\theta+1}}+
\sum_{l\in B-\lbrace 1\rbrace}{n_{F_l}\over d^{F_l}_{\theta}}+\sum_{F\in D_i^M}n_F}&\text{if $1\in A$} \\
\displaystyle{\sum_{l\in A-\lbrace 1\rbrace}{n_{F_l}\over
d^{F_l}_{\theta+1}}+\sum_{l\in B}{n_{F_l}\over
d^{F_l}_{\theta}}+(e^{F_1}_{\theta}-1){n_{F_1}\over
d^{F_1}_{\theta}}+\sum_{F\in D_i^M}n_F}&\text{if $1\in B$}
\end{cases}
\]
}

\noindent Let $(l_1,l2)\in A\times B$ and recall that
$\displaystyle {n_F\over d^{F}_{\theta+1}}$ (reps. $\displaystyle
{n_F\over d^{F}_{\theta}}$) does not depend on $F\in \cup_{l\in
A}P_{i_l}^{N_l}$ (resp. $F\in \cup_{l\in B}P_{i_l}^{N_l}$). In particular,
if we denote by $a$ (resp. $b$) the cardinality of $A$ (resp.
$B$), then we have:

\small{
\[
{\rm deg}_yQ_{f_y}(M,i)=
\begin{cases}

\displaystyle{a{n_{F_{l_1}}\over d^{F_{l_1}}_{\theta+1}}+(b-1){n_{F_{l_2}}\over d^{F_{l_2}}_{\theta}}+
\sum_{F\in D_i^M}n_F}&\text{if $1\in B$} \\
\displaystyle{(a-1){n_{F_{l_1}}\over d^{F_{l_1}}_{\theta+1}}+b {n_{F_{l_2}}\over
d^{F_{l_2}}_{\theta}}+(e^{F_1}_{\theta}-1)({n_{F_1}\over
d^{F_1}_{\theta}})+\sum_{F\in D_i^M}n_F}&\text{if $1\in A$}
\end{cases}
\]
}




\noindent Note also that if $B\not=\emptyset$ then $\displaystyle{{n_{F_{l_2}}\over d^{F_{l_2}}_{\theta}}={n_{F_1}\over d^{F_1}_{\theta}}}$,
on the other hand, if $B=\emptyset$, then $1\in A$.  In particular we get the following:

$$
{\rm deg}_yQ_{f_y}(M,i)=a{n_{F_{l_1}}\over d^{F_{l_1}}_{\theta+1}}+(b-1)
{n_{F_{1}}\over d^{F_{1}}_{\theta}}+\sum_{F\in D_i^M}n_F
$$

\noindent The above results can be stated as follows:

\begin{teorema}{\rm  Let $P_i^M$ be a point of  $T(f)$ and
    assume that $P_i^M\notin {\rm Top}(f)$. Let
$(P_{i_l}^{N_l}))_{1\leq l\leq t}$ be the set of points that
strictly dominate $P_i^M$ and let $\theta$ be the smallest integer
such that for all $F\in P_i^M, \displaystyle M\leq
{m^{F}_{\theta}\over n_F}$. Fix $F_l\in P_{i_l}^{N_{i_l}}$ for all $1\leq l\leq t$ and let $A$ (resp. $B$) be the set of $1\leq
l\leq t$ for which $M = \displaystyle{m^{F}_{\theta}\over n_F} $
(resp. $M <\displaystyle{m^{F}_{\theta}\over n_F}$) for all $F\in
\cup_{l\in A}P_{i_l}^{N_l}$ (resp. $F\in \cup_{l\in
B}P_{i_l}^{N_l}$). Let $(l_1,l_2)\in A\times B$. Let $F_{l_1}\in P_{i_{l_1}}^{N_{l_1}}$ and $F_{l_2}\in P_{i_{l_2}}^{N_{l_2}}$. If $a$ (resp.
$b$) denotes the cardinality of $A$ (resp. $B$) then the component
$Q_{f_y}(M,i)$ of $f_y$ satisfies the following:

\medskip

i) ${\rm deg}_yQ_{f_y}(M,i)=a{n_{F_{l_1}}\over d^{F_{l_1}}_{\theta+1}}+(b-1)
{n_{F_{l_2}}\over d^{F_{l_2}}_{\theta}}+\sum_{F\in D_i^M}n_F$, and $\sum_{F\in D_i^M}n_F$
 is given by the formula of Lemma 6.4., where if $F\in D_i^M$, then $h_F$ is
 either $\theta-1$ or $\theta$ depending on $M>\dfrac{m^F_{h_F}}{n_F}$ or
 $M=\dfrac{m^F_{h_F}}{n_F}$.

\medskip

ii) For all irreducible component $P$ of $Q_{f_y}(M,i)$ and for all
$F\in P_i^M, c(F,P)=M$.

\medskip

iii) For all irreducible component $P$ of $Q_{f_y}(M,i)$ and for all
$F\notin P_i^M$, c$(F,P)=c(F,P_{i}^M)<M$, where we recall that $c(F,P_i^M)$ is the
contact of
$F$ with any element
of $P_i^M$.

\medskip

iv) For all $1\leq k\leq \xi(f)$:

- If $f_k\in P_i^M$ then
int$(f_k,Q_{f_y}(M,i))=S(m^{f_k},M)\displaystyle{{{\rm
deg}_Q{f_y}(M,i)}\over {n_{f_k}}}$.

- If $f_k\notin P_i^M$ then
int$(f_k,Q_{f_y}(M,i))=S(m^{f_k},c(f_k,P_i^M))\displaystyle{{{\rm
deg}_Q{f_y}(M,i)}\over {n_{f_k}}}$, where $c(f_k,P_i^M)$ is the contact
of $f_k$ with any $F\in P_i^M$.$\blacksquare$

}
\end{teorema}

\noindent In the following we shall consider the case where $P_i^M$
is a top point of $T(f)$.

\begin{lema}{\rm Suppose that $P_i^M=D_i^M\in {\rm Top}(f)$, and let
    $F\in P_i^M$. We have the following:
$$
Q_{f_y}(M,i)={R_{=M}(F,f_y)}
$$
}
\end{lema}

\begin{demostracion}{.} By Lemma 7.2.,  deg$_yR_{>M}(G,f_y)=1$ for all
  $D_i^M$. Our assertion follows from Lemma 6.7.$\blacksquare$
\end{demostracion}

\noindent Let $P_i^M=D_i^M=\lbrace F_1,\ldots,F_r\rbrace$, and recall, by Proposition 2.4.,  that the 
sequence $(F_1,\ldots,F_r)$ is either
equivalent, or almost equivalent.

\begin{teorema}{\rm Let $P_i^M=\lbrace F_1,\ldots,F_r\rbrace \in{\rm Top}(f)$
and assume that $n_{F_1}={\rm max}_{1\leq k\leq
      r}n_{F_k}$. We have the following:

i) If $(F_1,\ldots,F_r)$ is equivalent with $M>\displaystyle{m^{F_1}_{h_{F_1}}\over n_{F_1}}$,
 then
deg$_yQ_{f_y}(M,i)= (r-1)n_{F_1}$.

ii) If  $(F_1,\ldots,F_r)$ is equivalent with $M=\displaystyle{m^{F_1}_{h_{F_1}}\over n_{F_1}}$, then
deg$_yQ_{f_y}(M,i)=
(r-1)n_{F_1}+(e^{F_1}_{h_{F_1}}-1)\displaystyle{n_{F_1}\over d^{F_1}_{h_{F_1}}}$.

iii) If $(F_1,\ldots,F_r)$ is
   almost equivalent, then
deg$_yQ_{f_y}(M,i)=(r-1)n_{F_1}$.
}
\end{teorema}
\begin{demostracion}{.} It follows from Lemma 7.8. that
deg$_yQ_{f_y}(M,i)={\rm deg}_y R_{=M}(F_1,f_y)$. Now the hypothesis of i) and ii) implies that
$n_{F_k}=n_{F_1}$ for all $k=2,\ldots,r$. Hence i) results from Proposition 5.2. and ii) results
from Proposition 5.3. Assume that $(F_1,\ldots,F_r)$ is almost equivalent, and that,
without loss of generality, $(F_1,F_3,\ldots,F_r)$ is equivalent. Since
$Q_{f_y}(M,i)=R_{=M}(F_1,f_y)=R_{=M}(F_2,f_y)$, then iii) results from Proposition 5.2.$\blacksquare$
\end{demostracion}

\begin{nota}{\rm When $P_i^M=D_i^M\in {\rm Top}(f)$, the numbers $a$ and $b$ of Theorem 7.7.
are zero. The reader may verify that the two formulas of Theorem 7.7. and Theorem 7.9.
coincide.$\blacksquare$
}
\end{nota}

\begin{exemple}{\rm i) Delgado's result: Let $f=f_1.f_2$. In [5],
in order to generalize Merle's Theorem, F. Delgado uses the arithmetic of the semi-group of $f$. His result is a
particular case of Theorem 7.7. More precisely, let $n_i={\rm
  deg}_yf_i,i=1,2$ and let $M=c(f_1,f_2)$, $I={\rm
int}(f_1,f_2)$. Let  $\theta$ be the smallest  integer such that $
M\leq {\dfrac{m^i_{\theta}}{n_{i}}}, i=1,2$. We have:

$$
f_y=(\prod_{k=1}^{\theta-1}Q_{f_y}(\dfrac{m^1_{k}}{n_{f_1}
},1)).\bar{f_y}
$$

\noindent where the properties $Q_{f_y}(\dfrac{m^1_{k}}{n(f_1)},1)$ are
given in the table 0), while those of the components of $\bar{f_y}$
are given in the tables 1), 2), 3), depending on the position of $M$
on $T(f)$. Note that $c(f_j,P)$ means the contact of $f_j$ with an
irreducible component of $Q_{f_y}(M,i)$.

 \bigskip
\noindent 0)


 \begin{tabular}{|*{4}{c|} l r | }
        \hline
    $Q$     & $Q_{f_y}({m^1_1\over {n_1}},1)$&...& $Q_{f_y}({m^1_{\theta-1}\over {n_1}},1)$ \\ \hline

${\rm deg}_yQ$ & $(e^1_1-1){{n_1}\over {d^{f_1}_1}}$&... &
$(e^1_{\theta-1}-1){{n_1}\over d^{f_1}_{\theta-1}}$ \\ \hline

$c(f_1,P), {\rm int}(f_1,Q)$ & ${m^{f_1}_1\over{n_1}},(e^1_1-1)r^1_1$&... &
${m^{f_1}_{\theta-1}\over {n_1}}, (e^1_{\theta-1}-1)r^1_{\theta-1}$ \\ \hline

$c(f_2,P), {\rm int}(f_2,Q)$ & ${m^{f_1}_1\over {n_1}},(e^1_1-1)r^2_1$&... &
${m^{f_1}_{\theta-1}\over {n_1}}, (e^1_{\theta-1}-1)r^2_{\theta-1}$ \\ \hline
\end{tabular}

\medskip

\noindent With the notations of  Theorem 7.7., for all $1\leq
i\leq \theta-1$, we have: $P_1^{\dfrac{m^1_i}{n_1}}=\lbrace f_1,f_2\rbrace, a=1, b=0$.


\bigskip

 \noindent 1) $M\not=\dfrac{m^i_{\theta}}{n(f_i)},i=1, 2$.

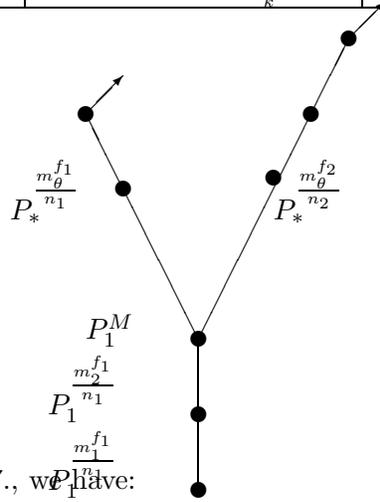
\begin{tabular}{|*{4}{c|} l r|}
        \hline
    $Q$ &    $Q_{f_y}(M,1)$ & $Q_{f_y}({m_k^{f_1}\over {n_1}},*),
\theta\leq k \leq h_{f_1}$& $Q_{f_y}({m_k^{f_2}\over {n_2}},*),\theta\leq k \leq h_{f_2}$\\ \hline

${\rm deg}_yQ$ &  ${{n_1}\over d^{f_1}_{\theta}}={{n_2}\over d^{f_2}_{\theta}}$ &
$(e^{f_1}_{k}-1){{n_1}\over {d^{f_1}_{k}}}$ & $(e^{f_2}_{k}-1){{n_2}\over {d^{f_2}_{k}}}$\\ \hline

$c(f_1,P), {\rm int}(f_1,Q)$ &$M, {{I}\over d^{f_1}_{\theta}}{{(n_1}\over {n_2}}$ &
${m^{f_1}_{k}\over {n_1}}, (e^{f_1}_{k}-1)r^{f_1}_{k}$ & $M, (e^{f_2}_{k}-1){{I}\over {d^{f_2}_{k}}}$ \\
\hline

$c(f_2,P), {\rm int}(f_2,Q)$ & $M, {{I}\over d^{f_2}_{\theta}}{{n_2}\over {n_1}}$ & $M,
(e^{f_1}_{k}-1){{I}\over {d^{f_1}_{k}}}$ & ${m^{f_2}_{k}\over {n_2}}, (e^{f_2}_{k}-1)r^{f_2}_{k}$ \\ \hline
\end{tabular}

\unitlength=1cm
\begin{center}
\begin{picture}(6,4)(3,-2)
\put(6,-4){\line(0,10){2}} 
\put(6,-4){\circle*{.2}}
\put(6,-3){\circle*{.2}}
\put(6,-2){\line(1,2){2}}  
\put(6,-2){\circle*{.2}}
\put(6,-2){\line(-1,2){1.5}} 
\put(5,0){\circle*{.2}}
 \put(4.5,1){\circle*{.2}}\put(4.5,1){\vector(1,1){0.5}}
 \put(7,0.15){\circle*{.2}}
\put(7.5,1){\circle*{.2}}
\put(8,2){\circle*{.2}}\put(8,2){\vector(1,1){0.5}}
\put(3.5,-0.4){\shortstack{$P_{*}^{{m^{f_1}_{\theta}\over n_1}}$}}
\put(7,-0.4){\shortstack{$P_{*}^{{m^{f_2}_{\theta}\over n_2}}$}}
\put(4.5,-2){\shortstack{$P_1^M$}}
\put(4,-4){\shortstack{$P_1^{{m^{f_1}_1}\over {n_1}}$}}
\put(4,-3){\shortstack{$P_1^{{m^{f_1}_2}\over{n_1}}$}}
\end{picture}
\end{center}

\bigskip
\bigskip
\bigskip
\bigskip

\noindent With the notations of Theorem 7.7., we have:

$P_1^M=\lbrace f_1,f_2\rbrace, A=\lbrace f_1\rbrace, B=\lbrace
f_2\rbrace, a=b=1$

$P_*^{{m_k^{f_1}}\over {n_1}}=\lbrace f_1\rbrace, \theta\leq k\leq
h_{f_1}: a=1,b=0$, $P_*^{{m_k^{f_2}}\over {n_2}}=\lbrace f_2\rbrace, \theta\leq k\leq h(f_2): a=1,b=0$

\noindent 2) $M=\dfrac{m^{f_1}_{\theta}}{n_1}< \dfrac{m^{f_2}_{\theta}}{n_2}$.

   \begin{tabular}{|*{4}{c|} l r|}
        \hline
    $Q$     & $Q_{f_y}(M,1)$ & $Q_{f_y}({m_k^{f_1}\over{n_1}},*), \theta+1\leq k \leq
h_{f_1}$ & $Q_{f_y}({m_k^{f_2}\over{n_2}},*),\theta\leq k \leq h_{f_2}$\\ \hline

${\rm deg}_yQ$ &
${{n_1}\over d^{f_1}_{\theta+1}}={{n_2}\over {d^{f_2}_{\theta}}}+(e^{f_1}_{\theta}-1){{n_1}\over
d^{f_1}_{\theta}}$&
$(e^{f_1}_{k}-1){{n_1}\over {d^{f_1}_{k}}}$ & $(e^2_{k}-1){{n_2}\over d^{f_2}_{k}}$\\ \hline

$c(f_1,P), {\rm int}(f_1,P_i)$ & $M,e^{f_1}_{\theta}r^{f_1}_{\theta}$ &
${m^{f_1}_k\over {n_1}}, (e^{f_1}_{k}-1)r^{f_1}_{k}$ & $M, (e^2_{k}-1){I\over d^{f_2}_{k}}$ \\
\hline

$c(f_2,P), {\rm int}(f_2,P_i)$ &$M,
{I\over d^{f_1}_{\theta+1}}=e^{f_1}_{\theta}r^{f_1}_{\theta}{{n_2}\over {n_1}}$ & $M,
(e^{f_1}_{k}-1){I\over d^{f_1}_{k}}$ & ${m^{f_2}_{k}\over {n_2}}, (e^{f_2}_{k}-1)r^{f_2}_{k}$ \\ \hline
\end{tabular}

\medskip

\noindent With the notations of Theorem 7.7., we have $P_1^M=\lbrace f_1,f_2\rbrace,
P_*^{m^{f_1}_k\over{n_1}}=\lbrace f_1\rbrace$ for all $\theta+1\leq k\leq h_{f_1}$, and
$P_*^{m^{f_2}_k\over{n_2}}=\lbrace f_2\rbrace$ for all $\theta\leq k\leq h_{f_2}$.

\bigskip
\noindent 3) $M=\dfrac{m^{f_1}_{\theta}}{n_1}=\dfrac{m^{f_2}_{\theta}}{n_2}$.

 \begin{tabular}{|*{4}{c|}l r|}
        \hline
    $Q$     & $Q_{f_y}(M,1)$ & $Q_{f_y}({m^{f_1}_k\over {n_1}},*), \theta+1\leq k \leq
h_{f_1}$ & $Q_{f_y}({m^{f_2}_k\over {n_2}},*), \theta+1\leq k \leq h_{f_2}$\\ \hline

${\rm deg}_yQ$ &
${{n_1}\over {d^{f_1}_{\theta+1}}}+(e^{f_1}_{\theta}-1){{n_1}\over d^{f_1}_{\theta}}$ &
$(e^{f_1}_{k}-1){{n_1}\over {d^{f_1}_{k}}}$ & $(e^{f_2}_{k}-1){{n_2}\over {d^{f_2}_{k}}}$\\ \hline

$c(f_1,P), {\rm int}(f_1,Q)$ & $M, (2e^{f_1}_{\theta}-1)r^{f_1}_{\theta}$
& ${m^{f_1}_{k}\over {n_1}}, (e^{f_1}_{k}-1)r^{f_1}_{k}$ & $M, (e^{f_2}_{k}-1){I\over d^{f_2}_{k}}$  \\
\hline

$c(f_2,P), {\rm int}(f_2,Q)$ & $M, (2e^{f_2}_{\theta}-1)r^{f_2}_{\theta}$
& $M, (e^{f_1}_{k}-1){I\over d^{f_1}_{k}}$ & ${m^{f_2}_{k}\over {n_2}}, (e^{f_2}_{k}-1)r^{f_2}_{k}$  \\
\hline
\end{tabular}

\medskip

\noindent With the notations of Theorem 7.7., we have $P_1^M=\lbrace f_1,f_2\rbrace,
P_*^{m^{f_1}_k\over{n_1}}=\lbrace f_1\rbrace$ for all $\theta+1\leq k\leq h_{f_1}$, and
$P_*^{m^{f_2}_k\over{n_2}}=\lbrace f_2\rbrace$ for all $\theta+1\leq k\leq h_{f_2}$.

\bigskip

\begin{exemple} {\rm i) $f=f_1.f_2$ and $f_1=(y^2-x^3)^2-x^5y,
f_2=(y^2-x^3)^2+x^5y$. We have $n_{f_1}=n_{f_2}=n=4,\underline{r}^{f_1}=\underline{r}^{f_2}=\underline{r}=(4,6,13),
\underline{d}^{f_1}=\underline{d}^{f_2}=\underline{d}=(4,2,1),\underline{m}^{f_1}=
\underline{m}^{f_2}=\underline{m}=(4,6,7)$, and $c(f_1,f_2)={\dfrac{7}{4}}$. The tree model of $f$ is
given by:

 \unitlength=1cm
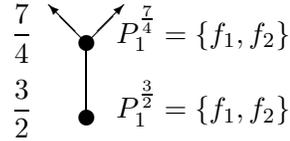
\begin{center}
\begin{picture}(6,2)(3,0.8)
\put(6,1){\line(0,10){1}} 
\put(6,1){\circle*{.2}}
\put(6,2){\circle*{.2}}
\put(6,2){\vector(-1,1){0.5}}
\put(6,2){\vector(1,1){0.5}}
\put(5,1){\shortstack{${\dfrac {3}{2}}$}} \put(5,2){\shortstack{$\dfrac{7}{4}$}}
\put(6.4,1){\shortstack{$P_1^{3\over 2}=\lbrace f_1,f_2\rbrace$}}
\put(6.4,2){\shortstack{$P_1^{7\over 4}=\lbrace f_1,f_2\rbrace$}}
\end{picture}
\end{center}

\noindent Note that $X({\dfrac{3}{2}},1)=\lbrace f_1\rbrace$, $X({\dfrac{3}{2}},2)=\lbrace
f_2\rbrace$, and $X({\dfrac{7}{4}},1)=P_1^{7\over 4}=\lbrace f_1,f_2\rbrace$. In particular  $c({\dfrac{3}{2}},1)=1$
and $c({\dfrac{7}{4}},1)=2$. With the notations of  Theorem 7.7.,
$f_y=Q({\dfrac{3}{2}},1)Q({\dfrac{7}{4}},1)=Q_1Q_2$, where:

\medskip

${\rm deg}_yQ_1={\dfrac{n}{d_2}}-{\dfrac{n}{d_1}}=1$ ($a=1,b=0$)

${\rm deg}_yQ_2={\dfrac{n}{d_3}}-{\dfrac{n}{d_2}}+n=2.4-2=6$ ($a=1, b=0$).

\medskip

\noindent Furthermore, for all irreducible component $P$ of $Q_1$
(resp. $Q_2$), $c(f_1,P)=c(f_2,P)={\dfrac{3}{2}}$ (resp.
$c(f_1,P)=c(f_2,P)={\dfrac{7}{4}}$). Finally,
int$(f_1,Q_1)=(e_1-1)r_1=r_1=6={\rm int}(f_2,Q_1)$ and
int$(f_1,Q_2)={\rm int}(f_1,f_2)+(e_2-1)r_2=39={\rm int}(f_2,Q_2)$.

\medskip

ii)  $f=f_1.f_2.f_3.f_4$ and $f_1=(y^2-x^3)^2-x^5y,
f_2=(y^2-x^3)^2+x^5y, f_3=(y^2-x^3)^2+x^5y-x^7$, and
$f_4=(y^2+x^3)^2-x^5y$: $\Gamma(f_i)=<4,6,13>=<n,r_1,r_2>,
i=1,2,3,4, m_1=6, m_2=7$, and $c(f_1,f_2)=c(f_1,f_3)=\dfrac{7}{4},
c(f_1,f_4)=\dfrac{3}{2}, c(f_2,f_3)=\dfrac{9}{4}, c(f_2,f_4)=\dfrac{3}{2}, c(f_3,f_4)=\dfrac{3}{2}$. The
tree model of $f$ is given by:

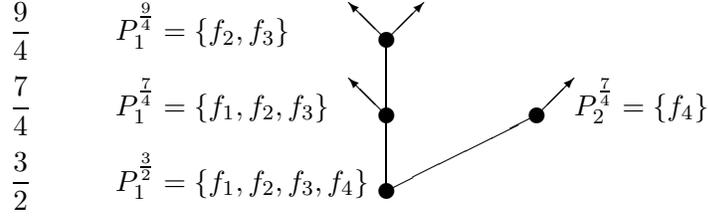
\unitlength=1cm
\begin{center}
\begin{picture}(6,4.5)(3,0.8)
\put(6,3){\circle*{.2}}
\put(6,3){\vector(1,1){0.5}}
\put(6,3){\vector(-1,1){0.5}}
\put(6,2){\circle*{.2}}
\put(6,2){\vector(-1,1){0.5}}
\put(6,1){\circle*{.2}}
\put(8,2){\circle*{.2}}
\put(8,2){\vector(1,1){0.5}}
\put(6,1){\line(0,10){2}} 
\put(6,1){\line(2,1){2}}
\put(1,1){\shortstack{$\dfrac{3}{2}$}}
\put(1,2){\shortstack{$\dfrac{7}{4}$}}
\put(1,3){\shortstack{$\dfrac{9}{4}$}}
\put(2.4,1){\shortstack{$P_1^{3\over 2}=\lbrace f_1,f_2,f_3,f_4\rbrace$}}
\put(2.4,2){\shortstack{$P_1^{7\over 4}=\lbrace f_1,f_2,f_3\rbrace$}}
\put(8.5,2){\shortstack{$P_2^{7\over 4}=\lbrace f_4\rbrace$}}
\put(2.4,3){\shortstack{$P_1^{9\over 4}=\lbrace f_2,f_3\rbrace$}}
\end{picture}
\end{center}

\bigskip

\noindent Note that $X_i(\dfrac{3}{2},1)=\lbrace f_i\rbrace, i=1,\ldots,4$,
 $X_1(\dfrac{7}{4},1)=\lbrace f_1,f_2\rbrace,
X_2(\dfrac{7}{4},1)=\lbrace f_1,f_3\rbrace$, $D_1^{7\over 4}=\lbrace
f_1\rbrace$, $X_2(\dfrac{7}{4},2)=P_2^{{7\over 4}}$ and $X_1(\dfrac{9}{4},1)=\lbrace
f_2,f_3\rbrace$. In particular, $c(\dfrac{3}{2},1)=1,
c(\dfrac{7}{4},1)=2=c(\dfrac{9}{4},1),c(\dfrac{7}{4},2)=1,$. Theorem 7.7.  implies that
$f_y=Q(\dfrac{3}{2},1)Q(\dfrac{7}{4},1)Q(\dfrac{7}{4},2)Q(\dfrac{9}{4},1)=Q_1Q_2Q_3Q_4$ with the
following properties:

\bigskip
  \begin{center}
     \begin{tabular}{|*{5}{c|}l r|}
        \hline
    $Q_i,{\rm deg}_yQ_i$     & $Q_1,3$& $Q_2,6$ & $Q_3,2$ & $Q_4,4$\\ \hline
   $c(f_1,P),{\rm int}(f_1,Q_i)$ & ${3\over 2}, 18$ & ${7\over 4},39$ &
${3\over 2},12$& ${7\over 4}, 26$ \\ \hline

 $c(f_2,P),{\rm int}(f_2,Q_i)$ & ${3\over 2}, 18$ & ${7\over 4},39$ &
${3\over 2},12$& ${9\over 4}, 28$ \\ \hline

$c(f_3,P),{\rm int}(f_3,Q_i)$ & ${3\over 2}, 18$ & ${7\over 4},39$ &
${3\over 2},12$&
${9\over 4}, 28$ \\ \hline

$c(f_4,P),{\rm int}(f_4,Q_i)$ & ${3\over 2}, 18$ & ${3\over 2}, 36$ &
${7\over 4}, 13$&
${3\over 2}, 24$ \\ \hline
\end{tabular}
  \end{center}

\noindent Where $c(F,P)$ means the contact of $F$ with an
irreducible component $P$ of $Q_i$.
}\end{exemple}

iii) Let $f=f_1.f_2.f_3$, where $f_1=(y^2-x^3)^2-x^5y, f_2=y^2-x^3$
and $f_3=y^2+x^3$. We have  $c(f_1,f_2)=\displaystyle{7\over 4},
c(f_1,f_3)=\displaystyle{3\over 2}=c(f_2,f_3),  {\rm int}(f_1,f_2)=13, {\rm int}(f_1,f_3)=12$ and
${\rm int}(f_2,f_3)=6$. The tree model of $f$ is given by:

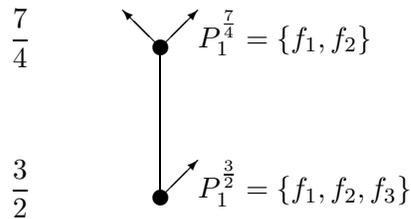
\unitlength=1cm
\begin{center}
\begin{picture}(6,4.5)(3,-1)
\put(6,2){\circle*{.2}}
\put(6,2){\vector(1,1){0.5}}
\put(6,2){\vector(-1,1){0.5}}
\put(6,0){\circle*{.2}}
\put(6,0){\vector(1,1){0.5}}
\put(6,0){\line(0,10){2}} 
\put(4,0){\shortstack{${\dfrac{3}{2}}$}}
\put(4,2){\shortstack{${\dfrac{7}{4}}$}}
\put(6.5,0){\shortstack{$P_1^{3\over 2}=\lbrace f_1,f_2,f_3\rbrace$}}
\put(6.5,2){\shortstack{$P_1^{7\over 4}=\lbrace f_1,f_2\rbrace$}}
\end{picture}
\end{center}
\noindent With the notations of Theorem 7.7., we have:

$ X({3\over 2},1)=\lbrace f_1,f_2\rbrace,
X({3\over 2},2)=\lbrace f_2,f_3\rbrace, D_1^{3\over 2}=\lbrace f_3\rbrace, c({3\over 2},1)=2$.

$ X({7\over 4},1)=\lbrace f_1,f_2\rbrace,
  c({7\over 4},1)=2$.

\noindent This gives us the following description:

\bigskip
\begin{center}

 \begin{tabular}{ |*{3}{c|}l r|}
 \hline $Q, {\rm deg}_yQ$ & $Q_{f_y}({3\over 2},1),5$ & $Q_{f_y}({7\over 4},1),2$ \\
\hline

$c(f_1,P), {\rm int}(f_1,Q)$ & ${7\over 4}, 26$ & ${3\over 2}, 18$\\ \hline

$c(f_2,P), {\rm int}(f_2,Q)$ & ${7\over 4}, 13$ & ${3\over 2},9$\\ \hline

$c(f_3,P), {\rm int}(f_3,Q)$ & ${3\over 2}, 12$ & ${3\over 2},9$\\ \hline
\end{tabular}
\end{center}

\noindent Where $c(F,P)$ means the contact of $F$ with an
irreducible component $P$ of $Q_i$.

iv) $f=f_1.f_2$, where $f_1=((y^2-x^3)^2-x^5y)^2+x^{10}(y^2-x^3)$
and $f_2=((y^2+x^3)^2-x^5y)^2+x^{22}(y^2+x^3)$. We have
$\Gamma(f_1)=<8,12,26,53>, \Gamma(f_2)=<8,12,26,57>, M=c(f_1,f_2)={3\over 2}$ and
$I={\rm int}(f_1,f_2)=96$. The tree model of $f$ is given by:

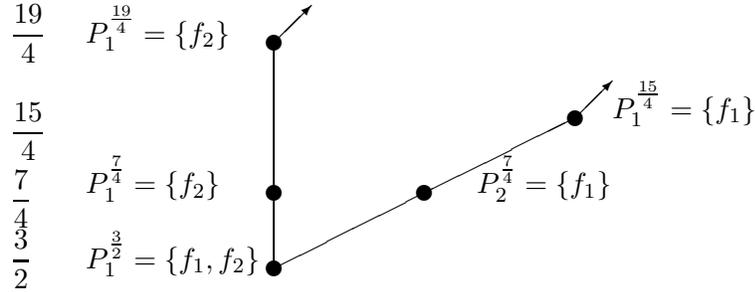
\unitlength=1cm

\begin{picture}(6,4.5)
\put(6,4){\circle*{.2}}
\put(6,2){\circle*{.2}}
\put(6,1){\circle*{.2}}
\put(10,3){\circle*{.2}}
\put(6,4){\vector(1,1){0.5}}
\put(8,2){\circle*{.2}}
\put(10,3){\vector(1,1){0.5}}
\put(6,1){\line(0,10){3}} 
\put(6,1){\line(2,1){4}}
\put(3.5,1){\shortstack{$P_1^{3\over 2}=\lbrace f_1,f_2\rbrace$}}
\put(2.5,1){\shortstack{${\dfrac{3}{2}}$}}
\put(3.5,2){\shortstack{$P_1^{7\over 4}=\lbrace f_2\rbrace$}}
\put(2.5,1.8){\shortstack{${\dfrac{7}{4}}$}}
 \put(10.5,3){\shortstack{$P_1^{15\over 4}=\lbrace f_1\rbrace$}}
\put(2.5,2.7){\shortstack{${\dfrac{15}{4}}$}}
\put(3.5,4){\shortstack{$P_1^{{19}\over 4}=\lbrace f_2\rbrace$}}
\put(2.5,4){\shortstack{${\dfrac{19}{4}}$}}
\put(8.7,2){\shortstack{$P_2^{7\over 4}=\lbrace f_1\rbrace$}}
\end{picture}

\noindent With the notations of Theorem 7.7., we have:

$X({\dfrac{3}{2}},1)=\lbrace f_1\rbrace,
X({\dfrac{3}{2}},2)=\lbrace f_2\rbrace, c({\dfrac{3}{2}}, 1)=1$.

$X({\dfrac{7}{4}},1)=\lbrace f_2\rbrace,
 c({\dfrac{7}{4}}, 1)=1, X({\dfrac{7}{4}},2)=\lbrace f_1\rbrace,
 c({\dfrac{7}{4}}, 2)=1$

$X({\dfrac{15}{4}},1)=\lbrace f_1\rbrace,
 c({\dfrac{15}{4}}, 1)=1, X({\dfrac{19}{4}},1)=\lbrace f_2\rbrace,
 c({\dfrac{19}{4}}, 1)=1$

\noindent This gives us the following description:

\bigskip

\begin{center}
  \begin{tabular}{ |*{6}{c|}l r|}
   \hline $Q, {\rm deg}_yQ$ & $Q_{f_y}({3\over 2},1),3$ & $Q_{f_y}({7\over 4},1),2$ & $Q_{f_y}({7\over 4},2)$,2 &
$Q_{f_y}({{15}\over 4},1),4$ & $Q_{f_y}({{19}\over 4},1),4$\\ \hline

$c(f_1,P), {\rm int}(f_1,Q)$ & ${3\over 2}, 36$ & ${3\over 2}, 24$ & ${7\over 4}, 26$ &
${{15}\over 4}, 53$ & ${3\over 2}, 48$\\ \hline

$c(f_2,P), {\rm int}(f_2,Q)$ & ${3\over 2}, 36$ & ${7\over 4}, 26$ & ${3\over 2}, 24$ &
${3\over 2}, 48$ & ${{19}\over 4}, 57$\\ \hline
\end{tabular}
  \end{center}

\noindent Where $c(F,P)$ means the contact of $F$ with an
irreducible component $P$ of $Q$.}\end{exemple}

\section{Factorization of the Jacobian}

\noindent Let $f=y^n+a_1(x)y^{n-1}+\ldots+a_n(x)$ and
$g=y^m+b_1(x)y^{m-1}+\ldots+b_m(x)$ be two monic reduced polynomials of
${\KK}((x))[y]$ and consider the Jacobian $J=J(f,g)$ of $f$ and
$g$. The aim of this Section is to give a factorization theorem of
$J$ in terms of the tree of $f.g$. Let to this end $T(f.g)$
be the tree of $f.g$ and let $f_1,\ldots,f_{\xi(f)}$ (resp.
$g_1,\ldots,g_{\xi(g)}$) be the irreducible components of $f$ (resp.
$g$) in ${\KK}((x))[y]$.

\begin{definicion}{\rm Let $P_i^M\in T(f.g)$.

i) We say that $P_i^M$ is an $f$-point if for all $1\leq k\leq \xi(g)$,
$g_k\notin P_i^M$ (equivalently $P_i^M$ is an $f$-point if for all $F\in
P_i^M$ and for all $1\leq k\leq \xi(g)$, $c(g_k,F) < M$).

ii) We say that $P_i^M$ is a $g$-point if for all  $1\leq k\leq
\xi(f)$, $f_k\notin P_i^M$ (equivalently $P_i^M$ is a $g$-point if
for all $F\in P_i^M$ and for all $1\leq k\leq \xi(f)$, $c(f_k,F) <
M$.

iii) We say that the point $P_i^M$ is
a mixed point if it is neither an $f$-point nor a $g$-point.

\vskip0.1cm

\noindent We
denote by $T_f$ (resp. $T_g$, resp. $T_m$) the set of $f$-points
(resp. $g$-points, resp. mixed points) of $T(f.g)$. Clearly
$T(f.g)=T_f\cup T_g\cup T_m$.}
\end{definicion}

\vskip0.5cm
 \unitlength=1cm
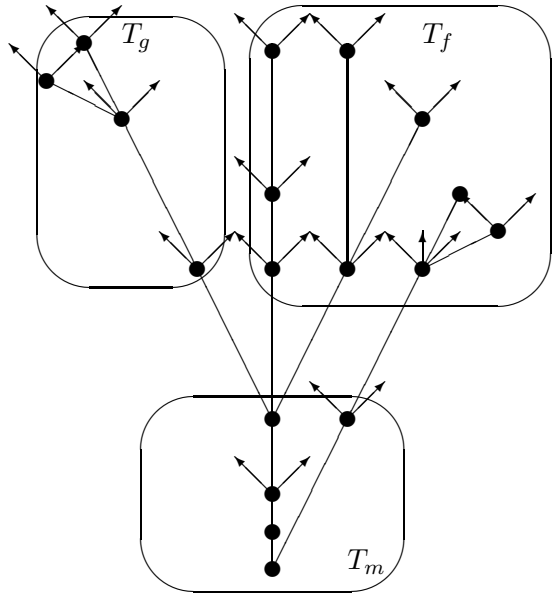
 \begin{center}
\begin{picture}(6,4)(3,-1)

\put(6,-4){\line(0,10){7}} 
\put(6,-4){\circle*{.2}}
\put(7,-4){\shortstack{$T_m$}}
 \put(6,-3.5){\circle*{.2}}
\put(6,0){\circle*{.2}}\put(6,0){\vector(1,1){0.5}}\put(6,0){\vector(-1,1){0.5}}
\put(6,1){\circle*{.2}}\put(6,1){\vector(1,1){0.5}}\put(6,1){\vector(-1,1){0.5}}
\put(6,2.9){\circle*{.2}}\put(6,2.9){\vector(1,1){0.5}}\put(6,2.9){\vector(-1,1){0.5}}
\put(6,-4){\line(1,2){2.5}}  
\put(7,-2){\circle*{.2}}\put(7,-2){\vector(1,1){0.5}}\put(7,-2){\vector(-1,1){0.5}}
\put(8.5,1){\circle*{.2}}
\put(6,-3){\circle*{.2}}\put(6,-3){\vector(1,1){0.5}}\put(6,-3){\vector(-1,1){0.5}}
 \put(6,-3){\oval(3.5,2.6)}
\put(6,-2){\line(1,2){2}}  
\put(8,2){\circle*{.2}}\put(8,2){\vector(1,1){0.5}}\put(8,2){\vector(-1,1){0.5}}
\put(8,3){\shortstack{$T_f$}}
\put(6,-2){\circle*{.2}}
\put(6,-2){\line(-1,2){2.5}} 
\put(5,0){\circle*{.2}}\put(5,0){\vector(1,1){0.5}}\put(5,0){\vector(-1,1){0.5}}
 \put(4,2){\circle*{.2}}\put(4,2){\vector(1,1){0.5}}\put(4,2){\vector(-1,1){0.5}}
\put(3.5,3){\circle*{.2}}\put(3.5,3){\vector(1,1){0.5}}\put(3.5,3){\vector(-1,1){0.5}}
\put(4,3){\shortstack{$T_g$}}
 \put(4,2){\line(-2,1){1}} 
\put(3,2.5){\circle*{.2}}\put(3,2.5){\vector(1,1){0.5}}\put(3,2.5){\vector(-1,1){0.5}}
\put(7,0){\line(0,3){3}} 
\put(7,0){\circle*{.2}}\put(7,0){\vector(1,1){0.5}}\put(7,0){\vector(-1,1){0.5}}
 \put(7.7,1.5){\oval(4,4)}
 \put(4.12,1.55){\oval(2.5,3.6)}
 \put(7,2.9){\circle*{.2}}\put(7,2.9){\vector(1,1){0.5}}\put(7,2.9){\vector(-1,1){0.5}}
\put(8,0){\line(2,1){1}}
\put(8,0){\circle*{.2}}\put(8,0){\vector(1,1){0.5}}\put(8,0){\vector(-1,1){0.5}}\put(8,0){\vector(0,1){0.5}}
\put(9,.5){\circle*{.2}}\put(9,0.5){\vector(1,1){0.5}}\put(9,0.5){\vector(-1,1){0.5}}
\end{picture}
\end{center}
\bigskip
\bigskip
\bigskip
\bigskip
\bigskip
\bigskip

\begin{lema}{\rm  Let $P_i^M,P_j^N\in
T(f.g)$, and assume that $P_j^N\geq P_i^M$.

i)  If $P_i^M\in T_f$ (resp. $P_i^M\in T_g$) then $P_j^N\in T_f$ (resp. $P_j^N\in
T_g)$.

ii)  if $P_j^N\in T_m$, then $P_i^M\in T_m$.}
\end{lema}

\begin{demostracion}{.} Easy exercise.$\blacksquare$
\end{demostracion}

\begin{lema}{\rm Let the notations be as above. If
    $T_f\not=\emptyset$ (resp. $T_g\not=\emptyset$), then
    Root$(J)\not=\emptyset$.
}

\end{lema}

\begin{demostracion}{.} Assume that $T_f\not=\emptyset$, and let $P=P_i^M\in T_f$.
Let $F\in P$ and let
$y_i(x),y_j(x)\in {\rm Root}(f)$ such that $c(y_i,y_j)=M$. By hypothesis,
  $M>{\rm max}_{z(x)\in{\rm Root}(g)}c(y_i,z)$. Now use Lemma 4.6.$\blacksquare$
\end{demostracion}

\noindent More generally, assume that $T_f\cup T_g\not=\emptyset$,
Propositions 5.3. and 5.4. and similar arguments as in Section 7.
led to the following factorization theorem of $J$.

\begin{teorema}{\rm  $J=\bar{J}.\prod_{P_i^M\in
      T_f}Q_{J}(M,i).\prod_{P_i^M\in T_g}Q_J(M,i)$, where for all $P_i^M\in T_f\cup
    T_g$, deg$_yQ_J(M,i)>1$. More precisely,  assume, without loss of generality, that $P=P_i^M\in T_f$ and let
$(P_{i_l}^{N_l}))_{1\leq l\leq t}$ be the set of points that strictly dominate
$P_i^M$. Let $\theta$ be the smallest integer such that
 $\displaystyle M\leq {m^{F}_{\theta}\over
n_F}$ for all $F\in P_i^M$. Let $A$ (resp. $B$) be the set of $1\leq l\leq t$ for
which $M = \displaystyle{m^{F}_{\theta}\over n_F} $ (resp. $M
<\displaystyle{m^{F}_{\theta}\over n_F}$) for all $F\in
\cup_{l\in A}P_{i_l}^{N_l}$ (resp. $F\in \cup_{l\in B}P_{i_l}^{N_l}$) and
let $(l_1,l_2)\in A\times B$. Let $F_{l_1}\in P_{i_{l_1}}^{N_{l_1}}$ and $F_{l_2}\in P_{i_{l_2}}^{N_{l_2}}$. If $a$ (resp. $b$) denotes the
cardinality of $A$ (resp. $B$) then the following hold:

\medskip

i) ${\rm deg}_yQ_{J}(M,i)=a.{n_{F_{l_1}}\over d^{F_{l_1}}_{\theta+1}}+(b-1).
{n_{F_{l_2}}\over d^{F_{l_2}}_{\theta}}+\sum_{F\in D_i^M}n_F$, and $\sum_{F\in D_i^M}n_F$  is given by the formula
of Lemma 5.6., where if $F\in D_i^M$, then $h_F$ is either $\theta-1$ or $\theta$ depending on $M>\dfrac{m^F_{h_F}}{n_F}$ or $M=\dfrac{m^F_{h_F}}{n_F}$.

\medskip

ii) For all irreducible component $P$ of $Q_{J}(M,i)$ and for all
$F\in P_i^M, c(F,P)=M$.

\medskip

iii) For all irreducible component $P$ of $Q_{J}(M,i)$ and for all
$F\notin P_i^M$ (this holds in particular when
$F=g_k, 1\leq k\leq \xi(g)$), c$(F,P)=c(F,P_{i}^M)<M$.

\medskip

iv) For all $1\leq k\leq \xi(f)$:

- If $f_k\in P_i^M$ then
int$(f_k,Q_{J}(M,i))=S(m^{f_k},M)\displaystyle{{{\rm
deg}_yQ_{J}(M,i)}\over {n_{f_k}}}$.

- If $f_k\notin P_i^M$ then
int$(f_k,Q_{J}(M,i))=S(m^{f_k},c(f_k,P_i^M))\displaystyle{{{\rm
deg}_yQ_{J}(M,i)}\over {n_{f_k}}}$.











}
\end{teorema}

\begin{demostracion}{.} The proof is similar to the proof of Theorem
  7.7.$\blacksquare$
\end{demostracion}

{\begin{corolario}{\rm Assume that $\xi(f)=1$, i.e. $f=f_1$ is an irreducible polynomial
of ${\KK}((x))[y]$, and let $M={\rm max}_{k=1}^{\xi(g)}c(f,g_k)$. Let $\theta$ be the
smallest integer such that $M< {m^f_{\theta}\over {n_f}}$. If $\theta <h_f$, then $J=J(f,g)=
\bar{J}.\prod_{k=\theta}^{h_f}J_k$, where for all $\theta\leq k\leq h_f$,

i) {\rm deg}$_yJ_k=(e^f_k-1)\displaystyle{n_f\over d^f_k}$.

ii) {\rm int}$(f,J_k)=(e^f_k-1)r^f_k$.

iii) For all irreducible component $P$ of $J_k$,
$c(f,P)=\displaystyle{m^f_k\over n_f}$.

iv) For all $1\leq j\leq \xi(g)$ and for all irreducible component $P$ of $J_k$,
$c(g_j,P)=c(g_j,f)$
}\end{corolario}

\bigskip

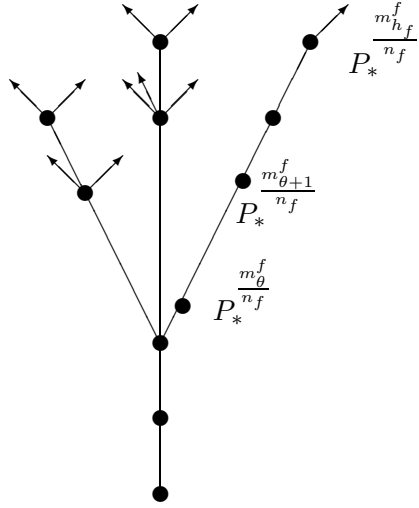
\unitlength=1cm
\begin{center}
\begin{picture}(6,4)(3,-2)
\put(6,-4){\line(0,10){6}} 
\put(6,-4){\circle*{.2}}
\put(6,1){\circle*{.2}}\put(6,1){\vector(1,1){0.5}}\put(6,1){\vector(-1,1){0.5}}
\put(6,1){\vector(-1,2){0.3}}
\put(6,2){\circle*{.2}}\put(6,2){\vector(1,1){0.5}}\put(6,2){\vector(-1,1){0.5}}
\put(6,-3){\circle*{.2}}
\put(6,-2){\line(1,2){2}}  
\put(6.3, -1.5){\circle*{.2}}
\put(6.7,-1.7){\shortstack{$P_{*}^{{m^{f}_{\theta}}\over{ n_f}}$}}

\put(6,-2){\circle*{.2}}
\put(6,-2){\line(-1,2){1.5}} 
\put(5,0){\circle*{.2}}\put(5,0){\vector(1,1){0.5}}\put(5,0){\vector(-1,1){0.5}}
 \put(4.5,1){\circle*{.2}}
\put(4.5,1){\vector(1,1){0.5}}\put(4.5,1){\vector(-1,1){0.5}}
 \put(7.1,0.15){\circle*{.2}}
\put(7.5,1){\circle*{.2}}
\put(8,2){\circle*{.2}}
\put(8,2){\vector(1,1){0.5}}
\put(8.5,1.6){\shortstack{$P_{*}^{{m^{f}_{h_f}}\over{ n_f}}$}}
\put(7,-0.4){\shortstack{$P_{*}^{{m^{f}_{\theta+1}}\over{ n_f}}$}}
\end{picture}
\end{center}
\bigskip
\bigskip
\bigskip
\bigskip
\begin{demostracion}{.} In fact, $T_f=\lbrace P_1^{m^f_{\theta}\over {n_f}},\ldots,
P_1^{m^f_{h_f}\over {n_f}}\rbrace$. The result  is consequently a particular case of Theorem
8.4.$\blacksquare$
\end{demostracion}

\section{Bad and good points on the tree of $f$}

\noindent Let $f=y^n+a_1(x)y^{n-1}+\ldots+a_n(x)$ be a monic reduced
polynomial of $\KK((x))[y]$, and let $f=f_1.\ldots.f_{\xi(f)}$
be the factorization of $f$ into irreducible components in
$\KK((x))[y]$. We shall assume that $f$ is generic in the following
sense: for all irreducible component $H$ of $f_y$, int$(f,H)\leq 0$.

\begin{definicion}{\rm Let $F,G$ be two monic polynomials of
    $\KK((x))[y]$, and let $H$ be an irreducible monic polynomial of
    $\KK((x))[y]$. We say that $H$ is regular (resp. irregular) with
    respect to $F$ if int$(F,H)\not=0$ (resp. int$(F,H)=0$). We define
    Reg$(G,F)$ (resp. Irreg$(G,F)$) to be the set of regular
    (resp. irregular) components of $G$ with
    respect to $F$. Let $\gamma(x)\in\KK((x^{1\over p})),
    p\in\mathbb{N}$. We say that $\gamma$ is regular (resp. irregular)
    with respect to $F$ if $O_xF(x,\gamma(x))\not= 0$
    (resp.$O_xF(x,\gamma(x))= 0$). If $G=F_y$, then we write
    Reg$(F)$
(resp. Irreg$(F)$) for Reg$(F_y,F)$ (resp. Irreg$(F_y,F))$.
}
\end{definicion}

\begin{lema}{\rm We have Irreg$(f,f_y)=\emptyset$.
}
\end{lema}

\begin{demostracion}{.} Let  $1\leq j\leq \xi(f)$ and let $y(x)\in{\rm
    Root}(f_j)$.  Let $M={\rm max}_{k\not=
    j}c(f_j,f_k)$. By Lemma 4.4., $M={\rm max}$ $c(f_j,H)$, where $H$ runs over
  the set of irreducible components of $f_y$. Since
  $f$ is generic, then $\sum_{y\not=\bar{y}\in{\rm
      Root}(f)}O_x(y-\bar{y})+M\leq 0$. If $M<0$, then
  $O_x(y-\bar{y})\leq M<0$ for all $\bar{y}\in{\rm
    Root}(f),\bar{y}\not= y$, in particular $\sum_{y\not=\bar{y}\in{\rm
      Root}(f)}O_x(y-\bar{y})<0$. If $M >0$, then $\sum_{y\not=\bar{y}\in{\rm
      Root}(f)}O_x(y-\bar{y})\leq -M <0$. Finally $O_xf_y(x,y(x))= \sum_{y\not=\bar{y}\in{\rm
      Root}(f)}O_x(y-\bar{y})<0$, in particular int$(f_j,f_y) <0$. This proves our assertion.$\blacksquare$
\end{demostracion}

\begin{definicion}{\rm Let $F,G$ be two monic polynomials of
    $\KK((x))[y]$, and let $H$ be an irreducible component of
    $G$. Assume that $H\in{\rm Irreg}(G,F)$ and let
$\gamma\in{\rm Root}(H)$. We have $F(x,\gamma(x))=\lambda+u(x)$ where
$\lambda\in\KK^*$ and $u(0)=0$. In particular, int$(F-\lambda,H)>0$,
hence $H\in{\rm Reg}(G,F-\lambda)$. We say that $\lambda$ is an
irregular value of $F$ with respect to $G$. We define irreg$(F,G)$ to
be
the set of irregular values of $F$ with respect to $G$. If $G=F_y$, then we write
    reg$(F)$
(resp. irreg$(F)$) for reg$(F_y,F)$ (resp. irreg$(F_y,F))$.

}
\end{definicion}

\begin{definicion}{\rm Let $P_i^M$ be a point of ${\rm Top}(f)$.

i) We say that
    $P_i^M$ is a good point if $H\in{\rm Reg}(f)$
    for some irreducible component $H$ of $Q_{f_y}(M,i)$.

ii) We say that
    $P_i^M$ is a bad point if $H\in{\rm Irreg}(f)$
    for some irreducible component $H$ of $Q_{f_y}(M,i)$.
}
\end{definicion}

\begin{lema}{\rm  Let $P_i^M$ be a point of ${\rm Top}(f)$.

i) If $P_i^M$ is a good point, then for all irreducible component
$H$ of $Q_{f_y}(M,i)$, $H\in{\rm Reg}(f)$.

ii) If $P_i^M$ is a bad point, then for all irreducible component
$H$ of $Q_{f_y}(M,i)$, $H\in{\rm Irreg}(f)$.

}
\end{lema}
\begin{demostracion}{.} i) By hypothesis, there is an irreducible component $\bar{H}$
of  $Q_{f_y}(M,i)$ such that int$(f,\bar{H})<0$. Let $H$ be an irreducible component of
$Q_{f_y}(M,i)$, and let $\gamma(x)$ (resp. $\bar{\gamma}(x)$) be a root of $H$
(resp. $\bar{H}$) such that max$_{i=1}^nc(\gamma,y_i)=M={\rm max}_{i=1}^nc(\bar{\gamma},y_i)$.
We have:
$$
 O_xf(x,\gamma(x))=\sum_{i=1}^nc(\gamma(x),y_i(x))=\sum_{i=1}^nc(\bar{\gamma}(x),y_i(x))=O_xf(x,\bar{\gamma}(x))
$$

\noindent in particular ${\rm int}(f,H)=\displaystyle{1\over n_H}O_xf(x,\gamma(x))<0$.

ii) The proof is similar to the proof of i).$\blacksquare$
\end{demostracion}

\section{Irregular values of a meromorphic curve}

\noindent Let the notations be as in Section 9, and let $P_i^M=\lbrace F_1,\ldots,F_r\rbrace$ be a
    bad point of ${\rm Top}(f)$. For all irreducible component $H$
    of $Q_{f_y}(M,i)$, int$(f,H)=0$, in particular, if $\gamma(x)\in
    {\rm Root}(H)$, then $f(x,\gamma(x))=\lambda+u(x)$, where
    $\lambda\in \KK^*$ and $u(0)=0$. In particular, $\lambda\in{\rm irreg}(f)$. Let
$\lbrace \lambda_1(M,i),\ldots,\lambda_{p(M,i)}(M,i)\rbrace$ be the set of
irregular values of $f$ obtained from the components of $Q_{f_y}(M,i)$
as above -more precisely $\lbrace \lambda_1(M,i),\ldots,\lambda_{p(M,i)}(M,i)\rbrace=
\lbrace {\rm inco}(f(x,\gamma(x)))| \gamma(x)\in{\rm Root}(Q_{f_y}(M,i))\rbrace$. We have the following:

\begin{proposicion}{\rm Assume that $n_{F_1}={\rm
      max}_{1\leq i\leq r}n_{F_i}$.

i) If $(F_1,\ldots,F_r)$ is equivalent and
$M>\displaystyle{m^{F_1}_{h_{F_1}}\over n_{F_1}}$, then $p(M,i)\leq
r-1$.

ii) If  $(F_1,\ldots,F_r)$ is equivalent and
$M=\displaystyle{m^{F_1}_{h_{F_1}}\over n_{F_1}}$, then $p(M,i)\leq
r$.

iii) If  $(F_1,\ldots,F_r)$ is almost equivalent,  then $p(M,i)\leq
r-1$.
}
\end{proposicion}

\begin{demostracion}{.} i) Let $H$ be an irreducible component of
  $Q_{f_y}(M,i)$. Since $c(H,F_1)=M >
  \displaystyle{m^{F_1}_{h_{F_1}}\over n_{F_1}}$, then $n_{F_1}$
  divides $n_H$. On the other hand, by Theorem 7.9., deg$_yQ_{f_y}(M,i)=(r-1)n_{F_1}$, In particular, $\xi(Q_{f_y}(M,i))\leq r-1$. This
  proves our assertion.

ii)  Let $H$ be an irreducible component of
  $Q_{f_y}(M,i)$. Since $c(H,F_1)=
  \displaystyle{m^{F_1}_{h_{F_1}}\over n_{F_1}}$, then
  $\displaystyle{n_{F_1}\over d^{F_1}_{h_{F_1}}}$
  divides $n_H$. More precisely, let $\gamma(x)=\sum_{p}c_px^{p\over
    n_H}\in{\rm Root}(H)$, then one of the following holds:

- The coefficient of $x^M$ in $\gamma(x)$  is nonzero, hence  $n_{F_1}$ divides
$n_H$. In this case, we say that $H$ is of type I.

- The coefficient of $x^M$ in $\gamma(x)$ is zero, then we say that $H$ is of type II.

\noindent Let $H_1,H_2$ be two irreducible components of type II of
$Q_{f_y}(M,i)$. If $\gamma_1(x)\in{\rm Root}(H_1)$
(resp. $\gamma_2(x)\in{\rm Root}(H_2)$), then
$c(y_i,\gamma_1)=c(y_i,\gamma_2)$, and inco$(y_i-\gamma_1)={\rm inco}(y_i-\gamma_2)$ for
all $y_i\in{\rm Root}(f)$. In particular, $H_1$ and $H_2$ give rise to
the same irregular value of $f$. On the other hand, by Theorem 7.9., deg$_yQ_{f_y}(M,i)=(r-1)n_{F_1}+(e^{F_1}_{h_{F_1}}-1)\displaystyle{{n_{F_1}\over d^{F_1}_{h_{F_1}}}}$, hence the number of
irreducible components of $Q_{f_y}(M,i)$ of type I is bounded by
$r-1$. This proves our assertion.

iii) The proof is similar to the proof of ii).$\blacksquare$

\end{demostracion}

\begin{corolario}{\rm Let $f$ be as above. The number of irregular
    values of $f$ is bounded by $\xi(f)$.
}
\end{corolario}

\begin{demostracion}{.} This results from Proposition 10.1.$\blacksquare$
\end{demostracion}





\begin{nota}{\rm Let the notations be above. If irreg$(f)$ has exactly $\xi(f)$
elements, then for all $P_i^M\notin {\rm Top}(f), D_i^M=\emptyset$.
More precisely, it follows from the proof of Proposition
  10.1. that the cardinality of irreg$(f)$ is bounded by

$$
\sum_{P_i^M\in{\rm Top}(f)}{\rm card}(P_i^M)
$$

\noindent In particular, if card(irreg($f$))$=\xi(f)$ then for
all $1\leq i\leq \xi(f)$, $f_i\in P_i^M$ for some bad point $P_i^M\in
{\rm Top}(f)$. Furthermore,
given a bad point $P_i^M=\lbrace
F_1,\ldots,F_r\rbrace\in{\rm Top}(f)$, the following holds:

i)  $(F_1,\ldots,F_r)$ is equivalent, and
$M=\displaystyle{M^{F_1}_{h_{F_1}}\over n_{F_1}}$.

ii) $Q_{f_y}(M,i)=H_1\ldots H_{r+1}$ and for all $i=1,\ldots,r$,
$H_i$ is irreducible of degree $n_{F_1}$ and $H_i$ is equivalent to
$F_1$. Furthermore, int$(H_{r+1},F_1\ldots
F_r)=(e^{F_1}_{h_{F_1}}-1)r^{F_1}_{h_{F_1}}$.

\noindent We do not have examples of meromorphic plane curve
satisfying the properties above, and we think that such an example
does not exist. More precisely, we think that the tree of a meromorphic
plane curve which is generic in its family must have at least one good point.
}
\end{nota}

\begin{nota}{\rm Suppose that $T(f)$ has only one bad point $P_i^M$,
    and that irreg$(f)$ has $\xi(f)$ elements (in
particular $P_i^M=\lbrace f_1,\ldots,f_{\xi(f)}\rbrace$). With the notations
of Remark 10.3., if $\xi(f)>1$ (resp. $\xi(f)=1$), then we have
int$(H_1,f)=0=\xi(f)r^{f_1}_{h_{f_1}}$
(resp. int$(H_1,f)=0=(e^{f_1}_{h_{f_1}}-1)r^{f_1}_{h_{f_1}}$), which
is a contradiction. This implies that  if $f$ has only one bad point, then
card(irreg($f$))$\leq\xi(f)-1$, and this bound is sharp (let
$f=y^4+x^{-1}y^2+y+1: \xi(f)=2$, $T(f)$ has one bad point and one good
point, and card(irreg($f$))$=1$). As a particular case, if $f$ is
irreducible, then irreg$(f)=\emptyset$. Note that if
$f\in\KK[x^{-1},y]$, then $f$ is irreducible in $\KK((x))[y]$ if and
only if $f(x^{-1},y)\in\KK[x,y]$ has one place at infinity. In this
case, the assertion above is a consequence of the Abhyankar-Moh theory.
}
\end{nota}





\noindent [1] S.S. Abhyankar.- On the semigroup of a meromorphic
curve, Part 1, in Proceedings, International Symposium on Algebraic
Geometry, Kyoto (1977), pp. 240-414.

\noindent [2] S.S. Abhyankar and A .Assi.- Jacobian of meromorphic curves,
Proc. Indian Acad. Sci. Math. Sci. 109, no. 2 (1999), 117-163.

\noindent [3] A. Assi.- Sur l'intersection des courbes m\'eromorphes,
C. R. Acad. Sci. Paris S\'er. I Math. 329 (1999), n$^0$ 7, 625-628.

\noindent [4] A. Assi.- Meromorphic plane curves,
Math. Z. 230 (1999), no. 1, 16-183.

\noindent [5] F. Delgado de la Mata.- A factorization theorem for
the polar of a curve with two branches, Compositio Mathematica, 92 (1994),
pp. 327-375.

\noindent [6] F. Delgado de la Mata.- An arithmetical  factorization for the critical point 
set of some map germs from ${\mathbb C}^2$ to ${\mathbb C}^2$, London Math. Soc.,
Lecture Note Ser.,  201 (1994),  pp. 61-100.

\noindent [7] E.R. Garc\'{i}a Barroso- Sur les courbes polaires d'une courbe plane r\'eduite, Proc. London Math. Soc. (3), 81, no.1 (2000)
pp. 1-28.

\noindent [8] T.C. Kuo and Y.C. Lu.- On analytic function germs of
two complex variables, Topology 16 (1977), pp. 299-310.

\noindent  [9] H. Maugendre.- Discriminant of a germ ${\mathbb C}^2\longmapsto {\mathbb C}^2$
and Seifert bref manifolds, J. London Math. Soc. (2), vol 59, no. 1 (999), pp. 207-226.

\noindent  [10] H. Maugendre.-Discriminant d'un germe $(g,f):({\mathbb C^2},0)\longmapsto ({\mathbb C}^2,0)$ et quotients de contact dans 
la rsolution de $fg$, Ann. Fac. Sci. Toulouse Math. (6) 7, no. 3 (1998), pp. 497-525.

\noindent [11] M. Merle.- Invariants polaires des courbes planes,
Inventiones Mathematics 41 (1977), pp. 299-310.

\noindent [12] C.T.C. Wall.- Singular points of plane curves, LMS Student Texts, 63, 
Cambridge University Press, 2004.

\noindent [13] O. Zariski.- Le probleme des modules pour les
branches planes, Lectures at Centre de Math\'ematiques, Ecole
Polytechnique, Notes by F. Kmety and M. Merle, 1973.

\end{document}